\documentstyle[11pt]{article}

\newtheorem{theorem}{Theorem}[section]

\newtheorem{assumption}[theorem]{Assumption}

\begin{document}
\title{Inverse problem for a parabolic system with two components by 
measurements of one component}
\author{Assia Benabdallah, \\
Laboratoire d'Analyse Topologie Probabilit\'{e}s/CMI, \\
U.M.R 6632, Universit\'{e} d'Aix-Marseille, Marseille, France\\
e-mail: assia@cmi.univ-mrs.fr,\\
Michel Cristofol, \\
Laboratoire d'Analyse Topologie Probabilit\'{e}s,\\
CNRS UMR 6632, Universit\'{e} d'Aix-Marseille, Marseille, France\\
e-mail:Michel.Cristofol@cmi.univ-mrs.fr,\\
Patricia Gaitan, \\
Laboratoire d'Analyse Topologie Probabilit\'{e}s,\\
CNRS UMR 6632, Universit\'{e} d'Aix-Marseille, Marseille, France\\
e-mail:Patricia.Gaitan@cmi.univ-mrs.fr,\\
Masahiro Yamamoto, \\
Department of Mathematics, University of Tokyo, \\
Komaba, Meguro, Tokyo 153, Japan \\
e-mail:  myama@ms.u-tokyo.ac.jp}

 \maketitle
\begin{abstract}
We consider a $2\times 2$ system of parabolic equations with first
and zeroth coupling and establish a Carleman estimate by extra data of 
only one component without data of initial values.  
Then we apply the Carleman estimate to inverse 
problems of determining some or all of the coefficients by 
observations in an arbitrary subdomain over a time interval 
of only one component and data of two components at a fixed positive 
time $\theta$ over the whole spatial domain.  The main results are 
Lipschitz stability estimates for the inverse problems.
For the Lipschitz stability, we have to assume some non-degeneracy
condition at $\theta$ for the two components and for it, we 
can approximately control the two components of the $2 \times 2$ 
system by inputs to 
only one component.  Such approximate controllability is proved also 
by our new Carleman estimate.
Finally we establish a Carleman estimate for a $3\times 3$ system
for parabolic equations with coupling of zeroth-order terms by one 
component to show the corresponding approximate 
controllability with a control to one component.

\end{abstract}
\date{}
%

\setcounter{equation}{0}
\section{Introduction and notations}
\label{sec: intro}

%
This article is devoted to the question of the identification of coefficients  for a reaction diffusion convection system of two equations in a bounded
domain, with the main particularity that we observe {\textit{only one component of the system}}.
Let $\Omega \subset {\bf R}^n$ be a bounded connected open set
with $C^2$-boundary $\partial\Omega$, and we set
$x = (x_1,..., x_n) \in {\bf R}^n$, $\partial_j = \frac{\partial}{\partial x_j}$,
$1 \le j \le n$, $\partial_t = \frac{\partial}{\partial t}$, $\nabla
= (\partial_1, ..., \partial_n)$, $\Delta = \sum_{j=1}^n \partial_j^2$.
For any fixed $T>0$, we set $\Omega_T= \Omega \times (0,T)$, 
$\Sigma_T=\partial \Omega \times (0,T)$ and we consider the
following $2\times 2$ reaction-diffusion-convection system~: 
\begin{equation}
\left\{ 
\begin{array}{lll}
\label{syst-UV}
\partial _{t}U =\Delta
U +a U+b V +A \cdot \nabla U +B \cdot \nabla V+f 
& \mbox{in} & \Omega _T, \\ 
\partial _{t}V =\Delta V +c U +d V +C \cdot \nabla
U +D \cdot \nabla V +g & \mbox{in} & \Omega _T, \\ 
U =h_1,\;V =h_2 & \mbox{on} & \Sigma _{T},\\
U(\cdot,0)=U_0,\, V(\cdot,0)=V_0& \mbox{in} & \Omega ,
\end{array}
\right.
\end{equation}
where $a,b,c,d$ are scalar functions and $A,B,C,D$ vectorial fields both 
defined on $\Omega$.
The boundary condition $h_i$ as well as $f,g$ shall be kept fixed. 
If we change
the reaction coefficients $b,c$ into $\widetilde{b},\widetilde{c}$, 
we let $(\widetilde{U},\widetilde{V})$ be the
solution of (\ref{syst-UV}) associated to $\widetilde{b},\widetilde{c}$ 
and 
$(\widetilde{U_0},\widetilde{V_0})$ for the
initial condition.  Let $\omega \subset \Omega$ be a non-empty
subdomain and $T>0$.  
We assume that we can measure both   
$$
U|_{\omega\times (0,T) } \mbox{ and }\quad 
(U,V)|_{\Omega\times \{\theta \} } .
$$   
at a time $\theta\in (0,T)$.

We set $\omega_T=\omega \times (0,T)$.  
For $m \in {\bf N}$, $1\leq p\leq \infty$, by $W^{m,p}(\Omega)$ and
$L^p(0,T;X)$ we denote the classical Sobolev space with the norm
$\Vert \cdot\Vert_{W^{m,p}(\Omega)}$, and
the space of $X$-valued $p$-Bochner integrable functions respectively
(e.g., [1]).  As usual we write $W^{0,p}(\Omega) = L^p(\Omega)$
and $H^m(\Omega) = W^{m,2}(\Omega)$ for $m \in {\bf N}$. 
We define a Banach space
$$
W_2^{m,\frac{m}{2}}(\Omega_T)=\{u:\Omega \times (0,T)\rightarrow {\bf R}
;\,\partial_x^{\alpha}\partial_t^{\alpha_{n+1}}u\in L^2(\Omega_T),\, 
 \mbox{for}\ |\alpha| + 2\alpha_{n+1} \le m\},
$$
with the norm
$$
\Vert u\Vert_{W_2^{m,\frac{m}{2}}(\Omega_T)}
= \sum_{|\alpha| + 2\alpha_{n+1} \le m}
\Vert \partial_x^{\alpha}\partial_t^{\alpha_{n+1}}u\Vert_{L^2(\Omega_T)}.
$$
Here $\alpha = (\alpha_1, \dots, \alpha_n)$ is a multi-index,
$\vert \alpha\vert = \alpha_1 + \cdots + \alpha_n$, 
$\partial_x^{\alpha} = \partial_1^{\alpha_1}\cdots \partial_n^{\alpha_n}$,
and the differentiation is to be understood in the weak sense. 
Let $M$ be an arbitrary positive constant. We denote by $\nu$
the outward unit normal to $\Omega$ and by $\overline{B}_{X}(0,r)$ the closed ball of a metric space $X$ centered on $0$ of radius $r$.

We pose the following assumptions.
\begin{assumption} 

(a) $a,b,\widetilde{b},c,\widetilde{c},d \in 
\overline{B}_{L^{\infty}(\Omega)}(0,M)$,

(b) $A,B,C,D \in \overline{B}_{L^{\infty}(\Omega) ^n}(0,M)$,

(c) $\omega \subset \Omega $ satisfies $\partial \omega \cap \partial
\Omega =\gamma \ $ and $\left| \gamma \right| $ $\neq 0$, and
$\omega$ is of class $C^{2}$,
 
(d) $|B(x)\cdot \nu (x)|\neq 0,\qquad x\in \gamma $,

(e) $B\in C^{2}(\overline{\omega  })^{n}$, 
$A\in C^{1}(\overline{\omega })^{n}$
and $b\in C^{2}(\overline{\omega })$,

(f) $\vert \widetilde{U}(\cdot,\theta) \vert, \thinspace
\vert \widetilde{V}(\cdot,\theta) \vert > \delta_0
\quad \mbox{on} \quad \overline{\Omega_T}$ with some constant 
$\delta_0>0$,

(g) $\Vert \widetilde{U} \Vert_{ C(\overline{\Omega_T} )}, 
\; \Vert \widetilde{V} \Vert_{ C(\overline{\Omega_T} )}\leq M$,
  
(h) $\Vert \widetilde{U} \Vert_{ C^3( \omega_T)}, \; \Vert \widetilde{V} \Vert_{ C^3( \omega_T)}\leq M$. 
\end{assumption}

If the functions and the coefficients appearing in (\ref{syst-UV}) 
satisfy sufficient smoothness and compatibility conditions, 
then Assumption 1.1 (g) and (h) are satisfied. 
By Ladyzenskaja, Solonnikov and Ural'ceva \cite{LSU} for example, 
we can describe such conditions, but
we are interested mainly in the inverse problem and we will not exploit
these conditions.    

Our first main result is the stability in determining the reaction 
coefficients $b,c$ :
\begin{theorem}\label{theo06}
Let $\theta\in (0,T)$ be fixed.  We suppose that Assumption 1.1 
is satisfied and that $(U,V)(\cdot,\theta) = (\widetilde{U},\widetilde{V})
(\cdot,\theta)$ in $\Omega$. 
Then there exists a constant $\kappa >0$  such that 
\begin{eqnarray}
\label{stability estimate}
\Vert b-\widetilde{b}\Vert_{L^2(\Omega)} + \Vert c-\widetilde{c}%
\Vert_{L^2(\Omega)} 
\leq \kappa \left( 
\Vert \partial_t(U - \widetilde{U}) \Vert_{ W_2^{2,1}(\omega_T )}
+ \Vert U  - \widetilde{U} \Vert_{ W_2^{2,1}(\omega_T)}\right)
\end{eqnarray}
\end{theorem}
The key ingredient to these stability results is a {\em global}
Carleman estimate for system~(\ref{syst-UV}).

Since the pioneer work of Bukhgeim-Klibanov \cite{BK:81}, 
Carleman estimates have been successfully used for the following 
problems:
\\
(i) the uniqueness and the stability in determining coefficients:
Especially for parabolic equations, see Benabdallah, Dermenjian and 
Le Rousseau \cite{BDL}, Benabdallah, Gaitan and Le Rousseau \cite{BGL}, 
Imanuvilov and Yamamoto \cite{IY:98}, \cite{IY}, 
Imanuvilov, Puel and Yamamoto [19],
Isakov \cite{I:98}, Klibanov \cite{K:92}, \cite{KT:04}
Klibanov and Timonov \cite{KT}, Yuan and Yamamoto \cite{YY:06} and the 
references therein.  For hyperbolic problems, among many works, 
we restrict ourselves to a few works such as Imanuvilov and Yamamoto 
\cite{IY:01}, Isakov \cite{Is1}, \cite{I:98},
Klibanov \cite{K:92}, Klibanov and Timonov \cite{KT}
and see the references also in Isakov \cite{I:98} and Klibanov and 
Timonov \cite{KT}.
\\
(ii) observability inequalities and related estimates: 
see Fursikov and Imanuvilov \cite {FI}, 
Imanuvilov \cite{I}, Isakov \cite{Is1}, \cite{I:98}, 
Kazemi and Klibanov \cite{KK}, 
Klibanov and Malinsky \cite{KM}.
Furthermore the exact controllability of linear systems is equivalent 
to the observability of the corresponding adjoint system and 
we can refer to \cite{FI}, \cite{I}.
Imanuvilov and Yamamoto \cite{IY} discuss the 
global exact zero controllability 
for a semilinear parabolic equation.  Also see 
Ammar-Khodja, Benabdallah and Dupaix \cite{ABD}, and
Ammar-Khodja, Benabdallah, Dupaix and Kostine \cite{ABDK1}, \cite{ABDK},
Gonz\'{a}lez-Burgos and  P\'{e}rez-Garc\'{i}a \cite{GPe2} 
for semilinear parabolic systems.

Apart from the last previous works quoted, the existing Carleman estimates require observations of all the 
components when we will discuss inverse problems for a system such as (\ref{syst-UV}).  It is very desirable to establish the stability for inverse
problems for a $2 \times 2$ parabolic system by means of only one 
component,
because for a reaction-diffusion system it may be frequently difficult to observe the both components.
There are not many papers devoted to such inverse problems for $2\times
2$ parabolic systems, and
we can refer, for instance, to  Cristofol, Gaitan and Ramoul \cite{CGR}. 

The article is organized as follows. In Section~\ref{sec: 1} we derive
a new Carleman estimate for system~(\ref{syst-UV}). In Section~\ref{sec: 3} we
prove the stability result. In Section~\ref{sec: 4} we will remove  
Assumption 1.1 (f) on positivity of $\tilde{U},\tilde{V}$ at a time $\theta >0$. Section~\ref{sec: 5} is devoted to some comments and open problems.
The appendices provide technical proofs of lemmata stated in Sections~\ref{sec: 1} and ~\ref{sec: 4}. We want to point that the Carleman estimate proved in Section~\ref{sec: 1} implies a {\textit{new approximate controllability 
result}} for a $2\times 2$ reaction-diffusion-convection system with 
{\textit{one localized control}}. As it will be seen in Section~\ref{sec: 5}, this result can be extended to a $3\times 3$ 
{\textit{reaction-diffusion system}}.



\section{Carleman estimate } \label{sec: 1}

%
\subsection{A Carleman estimate for a $2\times 2$ system 
by extra data of one compoment}
Let $(a_{ij})_{1\le i,j\le 2}\in
L^{\infty }(\Omega_T )$ and $(A_{ij})_{1\le i,j\le 2}\in 
L^{\infty }(\Omega_T)^{n}.$ Let 
$u_{0},v_{0}\in L^{2}(\Omega )$ and $f,g\in L^2(\Omega_T)$. 
Consider the following reaction-diffusion system with convection terms~: 
 \begin{equation}
\left\{ 
\begin{array}{lll}
\partial _{t}u=\Delta u+a_{11}\ u+a_{12}\ v+A_{11}\cdot\nabla u
+A_{12}\cdot\nabla v+f & \mbox{in} & 
\Omega _{T}, \\ 
\partial _{t}v=\Delta v+a_{21}\ u+a_{22}\ v+A_{21}\cdot\nabla u
+A_{22}\cdot\nabla v+g & \mbox{in} & 
\Omega _{T}, \\ 
u = v = 0 & \mbox{on} & \Sigma _{T}, \\ 
u(\cdot,0)=u_{0}, \thinspace v(\cdot,0)=v_{0} & \mbox{in} & \Omega.
\end{array}
\right.  \label{syst-uv}
\end{equation}
Uniqueness existence and stability results in solving an initial 
value-boundary value problem (\ref{syst-uv}) can be proved 
by the semigroup theory for example (e.g., \cite{LSU}, Pazy 
\cite{Pa}, Tanabe \cite{T}).
In particular it admits a unique solution $(u,v)\in 
C([0,T];L^2(\Omega ))^2\cap L^2(0,T;H^1_0(\Omega ))^2$. 

Our main interest is to derive a Carleman estimate of $(u,v)$
solution of (\ref{syst-uv}) by \textit{solely} observing $u$ in 
$\omega \times (0,T).$ 
We make the following main assumptions~:
\begin{assumption}
\label{ass2}

\noindent
(a) Let $\omega \subset \Omega $ with $\partial \omega \cap \partial
\Omega =\gamma \ $and $\left| \gamma \right| $ $\neq 0$.\\ 
(b) $|A_{12}(x,t)\cdot \nu (x)|\neq 0,\qquad (x,t)\in \gamma_T$, with $\gamma_T=\gamma \times (0,T)$,\\
(c) $\Vert A_{12}\Vert_{C^{2}(\overline{\omega_T})^{n}}$, 
$\Vert a_{12}\Vert_{C^{2}(\overline{\omega_T})}$,
$\Vert A_{11}\Vert_{C^{1}(\overline{\omega_T})^{n}} \le M$,
where $M>0$ is an arbitrarily fixed constant.
\end{assumption}

In the sequel $\kappa $ will denote a generic constant and their 
values may change from a line to others.  The dependence of $\kappa$ 
on $s$ will be specified. 
 
In this section, we prove: 
\begin{theorem}\label{newcarl}
Let $\tau \ge 1$ and $\omega \subset \Omega$ be a subdomain such
that $\overline{\omega} \subset \Omega$. 
Under Assumption 2.1, there exist $\alpha_{\omega} \in C^2
(\overline{\Omega})$ with $\alpha_{\omega} > 0$ on $\overline
{\Omega}$ and 
two positive constants $s_0$ and $\kappa$ which depend on $T, M, 
\Omega $, $\omega$, $\tau$ and the $L^{\infty}$-norms of 
$a_{ij},A_{ij}$, such that there exist positive constants 
$\kappa_1(s,\tau)$ and $\kappa$ such that the following Carleman
estimate holds 
$$
\int_{\Omega_T} (s\rho)^{\tau-1}e^{-2s\eta_{\omega}}
(\vert \partial_t u\vert^2 + \vert \partial_t v\vert^2
+ \vert \Delta u\vert^2 + \vert \Delta v\vert^2
+ (s\rho)^2\vert \nabla u\vert^2 
+ (s\rho)^2\vert \nabla v\vert^2
+ (s\rho)^4 \vert u\vert^2 + (s\rho)^4\vert v\vert^2)
$$
$$
\le \kappa_1(s,\tau)(\Vert u\Vert^2_{W_2^{2,1}(\omega_{T})}
+ \Vert f\Vert_{L^2(\omega_T)}^2) 
+ \kappa \int_{\Omega _{T}}(s\rho )^{\tau}e^{-2s\eta_{\omega}}
(|f|^{2}+|g|^{2})
$$
for all $s \ge s_0$ and any solution $(u,v)$ to (\ref{syst-uv}).
Here we set 
\begin{equation}
\eta_{\omega}(x,t) = \frac{\alpha_{\omega}(x)}{t(T-t)}, \quad
\rho(t) = \frac{1}{t(T-t)}.               \label{(2.2)}
\end{equation}
\end{theorem}

This is a Carleman estimate for a $2\times 2$ system with extra data 
in $\omega_T$ of only one component.
In \cite{ABD} and \cite{CGR}, it is assumed that $A_{11} = A_{12}
= 0$.  In that case, the proof can be completed by 
directly substituting $v$ by means of $u$ in $\omega_T$.
By the first-order coupling, we extra need Assumption 2.1 (a) and (b).
 
{\bf Proof of Theorem 2.2}
First we prove
\\
{\bf Lemma 2.3}
{\it Let $\omega \subset \Omega$ be a subdomain and 
$\partial\omega \cap \partial\Omega = \gamma$. We consider 
\begin{equation}
\sum^n_{j=1} p_j(x,t)\partial_j u(x,t) + q(x,t)u(x,t) = f(x,t), \quad
x\in \omega \subset\Omega, \thinspace 0 < t < T. 
                                       \label{(2.3)}
\end{equation}
Here $p_j, q \in L^{\infty}(0,T;C^1(\overline{\Omega}))$ for 
$1 \le j \le n$.  We set $p = (p_1,..., p_n)$ and let 
$\nu(x) = (\nu_1(x), ..., \nu_n(x))$ be the unit outward normal vector 
to $\partial\omega$ at $x$.
We assume that
\begin{equation}
\vert p(x,t)\cdot \nu(x)\vert \ne 0, \qquad x \in \overline{\gamma}, 
\thinspace 0 \le t \le T.               \label{(2.4)}      
\end{equation}
Let $u = u(x,t)$ satisfy (\ref{(2.3)}) and $u\vert_{\gamma\times (0,T)} = 0$.
Then there exist a subdomain $\omega' \subset \omega$ and a constant
$\kappa>0$, which is dependent on $p$ and $q$ and independent of 
$f$, such that 
$$
\Vert u\Vert_{L^2(\omega'_T)} \le \kappa\Vert f\Vert_{L^2(\omega'_T)}.
$$
}

{\bf Proof of Lemma 2.3.}
We set $x = (x_1,..., x_n) = (x',x_n)$ and $y = (y_1, ..., y_n)
= (y', y_n)$.  Without loss of generality, we can assume that
$$
\omega = \{ (x',x_n);\thinspace 
h(x') < x_n < h_1(x'), \thinspace \vert x'\vert < \rho\}
$$
and $\gamma = \{ (x',x_n); \thinspace x_n = h_1(x'), \thinspace
\vert x'\vert <\rho\}$.
Here $\rho > 0$ is sufficiently small and $h, h_1
\in C^2(\{ \vert x'\vert\le \rho\})$ satisfy $h = h_1$ on 
$\{ \vert x'\vert = \rho \}$.
We change independent variables $y' = x'$ and $y_n = x_n - h(x')$.
Then $\omega$ is transformed to
$$
\widetilde{\omega} = \{ (y',y_n); \thinspace
0 < y_n < (h_1-h)(x'), \thinspace \vert y'\vert < \rho\}.
$$
Set $\tilde{u}(y,t) = u(x,t)$, $\tilde{p}(y,t) = p(x,t)$, 
$\tilde{q}(y,t) = q(x,t)$, $\tilde{f}(y,t) = f(x,t)$,
$\tilde{\Gamma}_1 = \{ (y',0); \thinspace \vert y'\vert < \rho\}$
and $\tilde{\Gamma}_2 = \{ (y',y_n); \thinspace
y_n = (h_1-h)(y'), \thinspace \vert y'\vert < \rho\}$.
Then $\partial\tilde{\omega} = \tilde{\Gamma}_1 \cup \tilde{\Gamma}_2$,
$$
(\tilde{P}\tilde{u})(y,t) = \sum^{n-1}_{j=1} \tilde{p}_j(y,t)
\frac{\partial\tilde{u}}{\partial y_j}(y,t)
$$
\begin{equation}
+ \tilde{r}(y,t)\frac{\partial \tilde{u}}{\partial y_n}(y,t)
+ \tilde{q}(y,t)\tilde{u}(y,t) = \tilde{f}(y,t), \quad y\in \tilde{\omega},
\thinspace 0 < t < T                                 \label{(2.5)}
\end{equation}
where
$$
\tilde{r}(y,t) = \tilde{p}_n(y,t) - \sum^{n-1}_{j=1} \tilde{p}_j(y,t)
\frac{\partial h}{\partial y_j}(y,t),
$$
and
\begin{equation}
\tilde{u}(y',y_n,t) = 0, \qquad y_n = (h_1-h)(y'), \thinspace
\vert y'\vert < \rho, \thinspace 0 < t <T.                  \label{(2.6)}
\end{equation}
Moerover
$\nu(x)$ is parallel to $(\partial_1h(x'), ...., \partial_{n-1}h(x'), -1)$ 
on $\{ (x',x_n); \thinspace x_n = h(x'), \thinspace \vert x'\vert 
< \rho \}$.
Therefore, in terms of  (\ref{(2.4)}), without loss of generality, we can assume
that there exists a constant $\delta > 0$ such that 
$\tilde{r}(y',0,t) > 2\delta$ for $\vert y'\vert < \rho$ and $0 < t < T$.
We choose $\rho > 0$ sufficiently small, so that 
\begin{equation}
\tilde{r}(y,t) > \delta, \qquad y\in \tilde{\omega}, \thinspace
0 < t < T.                                           \label{(2.7)}
\end{equation}
Let $\tilde{\nu}(y) = (\tilde{\nu}_1(y), ..., \tilde{\nu}_n(y))$ be the 
unit outward normal vector to $\partial\tilde{\omega}$ at $y$.
Then $\tilde{\nu}(y)$ is parallel to $(0, ..., 0, -1)$ for $y \in 
\tilde{\Gamma}_1$ and to 
$\left( - \frac{\partial(h_1-h)}{\partial y_1}(y'), ...,
- \frac{\partial(h_1-h)}{\partial y_{n-1}}(y'), 1\right)$ for
$y \in \tilde{\Gamma}_2$.  

Hence, by choosing $h_1, h$ such that $\Vert h_1-h\Vert_{C^1(\{ \vert y'\vert\le\rho\})}$ is sufficiently small if necessary, by (\ref{(2.7)}) we have
\begin{equation}
\mbox{$\widetilde{\Gamma}_1 \subset \left\{ y \in \partial\tilde{\omega}; 
\thinspace
\sum_{j=1}^{n-1} \tilde{p}_j(y,t)\tilde{\nu}_j(y) + \tilde{r}(y,t)
\tilde{\nu}_n(y) \le 0 \right\}$ }                   \label{(2.8)}
\end{equation}
and
\begin{equation}
\tilde{u}(\cdot,t) = 0 \quad\mbox{on $\tilde{\Gamma}_2$}, 
\quad \tilde{\Gamma}_2 \subset \left\{ y \in \partial\tilde{\omega}; \thinspace
\sum_{j=1}^{n-1} \tilde{p}_j(y,t)\tilde{\nu}_j(y) + \tilde{r}(y,t)
\tilde{\nu}_n(y) > 0 \right\}.                          \label{(2.9)}
\end{equation}
For the proof of the lemma, it suffices to prove a Carleman estimate
for (\ref{(2.3)}), whose proof is similar for example to Lemma 3.2 
in \cite{IY:05}.  We set 
$$
\tilde{P}_0\tilde{u} = \widetilde{P}\widetilde{u}
- \tilde{q}\tilde{u} = \sum_{j=1}^{n-1} \tilde{p}_j(y,t)
\frac{\partial\tilde{u}}{\partial y_j}(y,t) + \tilde{r}(y,t)\frac{\partial\tilde{u}}
{\partial y_n}(y,t),
$$
$w=w(\cdot,t) =\tilde{u}(\cdot,t) e^{sy_n}$ and $Qw = e^{sy_n}\tilde{P}_0
(e^{-sy_n}w)$.  Then
$$
Qw = \tilde{P}_0w - s\tilde{r}(y,t)w.
$$
We arbitrarily fix $t \in [0,T]$.
Hence by integration by parts and (\ref{(2.6)}) - (\ref{(2.9)}) we obtain
\begin{eqnarray*}
&&\int_{\tilde{\omega}} \vert \tilde{P}_0\tilde{u}\vert^2 e^{2sy_n} dy
= \int_{\tilde{\omega}} \vert Qw\vert^2 dy\\
=&& \Vert \tilde{P}_0w\Vert^2_{L^2(\tilde{\omega})} 
+ s^2\Vert \tilde{r}(\cdot,t)w\Vert^2_{L^2(\tilde{\omega})}
- 2s\int_{\tilde{\omega}}\left( 
\sum_{j=1}^{n-1} \tilde{p}_j\frac{\partial w}{\partial y_j}(y,t) 
+ \tilde{r}\frac{\partial w}{\partial y_n}\right) \tilde{r}w dy\\
\ge &&s^2 \int_{\tilde{\omega}} \tilde{r}^2w^2 dy
- s\int_{\tilde{\omega}}
\left(\sum_{j=1}^{n-1} \tilde{p}_j\tilde{r}\frac{\partial w^2}{\partial y_j} 
+ \tilde{r}^2\frac{\partial w^2}{\partial y_n}\right) dy\\
\ge&& s^2\delta\int_{\tilde{\omega}} w^2 dy
+ s\int_{\tilde{\omega}}\left( \sum_{j=1}^{n-1} 
\frac{\partial(\tilde{p}_j\tilde{r})}{\partial y_j} 
+ \frac{\partial\tilde{r}^2}{\partial y_n}\right)w^2 dy
- s \left( \int_{\tilde{\Gamma}_1} + \int_{\tilde{\Gamma}_2}\right)
\tilde{r}\left( 
\sum_{j=1}^{n-1} \tilde{p}_j\tilde{\nu}_j + \tilde{r}\tilde{\nu}_n
\right)w^2 dS\\
\ge&& s^2 \int_{\tilde{\omega}} \left( \delta - \frac{\kappa_1}{s}
\right) w^2 dy.
\end{eqnarray*}
Henceforth $\kappa_j > 0$ depends on $\max_{1\le j \le n}
\Vert p_j\Vert_{C^1(\overline{\Omega_T})}$ and $\omega$.
Hence we have
$$
s^2\int_{\tilde{\omega}} \vert \tilde{u}\vert^2e^{2sy_n} dy
\le \kappa_2\int_{\tilde{\omega}} \vert \tilde{P}_0\tilde{u}\vert^2
e^{2sy_n} dy
$$
for all large $s > 0$.  Since
\begin{eqnarray*}
&& \int_{\tilde{\omega}} \vert \tilde{P}_0\tilde{u}\vert^2 
e^{2sy_n} dy
\le 2\int_{\tilde{\omega}} \vert \tilde{P}\tilde{u}\vert^2 e^{2sy_n}dy
+ 2 \int_{\tilde{\omega}} \vert \tilde{q}\tilde{u}\vert^2 e^{2sy_n} dy\\
\le &&2\int_{\tilde{\omega}} \vert \tilde{f}\vert^2 e^{2sy_n}dy
+ 2\Vert q\Vert^2_{C(\overline{\Omega_T})}
\int_{\tilde{\omega}} \vert \tilde{u}\vert^2 e^{2sy_n} dy,
\end{eqnarray*}
by choosing $s$ large such that $\frac{s^2}{\kappa_2} - 2\Vert q\Vert^2
_{C(\overline{\Omega_T})} \ge \frac{s^2}{2}$, we have
$$
s^2\int_{\tilde{\omega}} \vert \tilde{u}\vert^2e^{2sy_n} dy
\le \kappa_3\int_{\tilde{\omega}} \vert \tilde{f}\vert^2
e^{2sy_n} dy
$$
for all large $s > 0$.  Since $1 \le e^{2sy_n} \le e^{2s\kappa_4}$
for $y \in \tilde{\omega}$ where $\kappa_4 = \Vert h_1-h
\Vert_{C(\{ \vert y'\vert\le \rho\})}$, for all large $s>0$, we fix
$s > 0$ large and we have
$\Vert \tilde{u}(\cdot,t)\Vert_{L^2(\tilde{\omega})} \le
\kappa_5\Vert \tilde{f}(\cdot,t)\Vert_{L^2(\tilde{\omega})}$.
By integrating over $t \in (0,T)$, the proof of Lemma 2.3 is completed.

By (\ref{syst-uv}), we have
$$
A_{14}\cdot \nabla v + a_{12}v =
\partial _{t}u-\Delta u+a_{11}\ u+A_{13}\ .\nabla u+f   
\mbox{ in } \omega _{T}                                
$$
and
$$
v=0 \mbox{ on }\partial\Omega \times (0,T).
$$
In terms of Assumption 2.1 (b), we apply Lemma 2.3 so that we can
choose a subdomain $\omega' \subset \Omega$ such that
\begin{equation}
\Vert v\Vert_{L^2(\omega_T')} \le \kappa\Vert u\Vert
_{W_2^{2,1}(\omega'_T)}+\Vert f\Vert_{L^2(\omega '_T)}.                            \label{(2.10)}
\end{equation}
By \cite{FI} and \cite{I}, for $\omega'$, there exist $\beta_{\omega'}
\in C^2(\overline{\Omega})$ with $\beta_{\omega'} > 0$ on $\overline
{\Omega}$ and two positive constants $s_0$ and $\kappa$,
which depend on $ T, \Omega $, $\omega'$, $\tau$ and $L^{\infty}$ norms 
of $a_{ij},A_{ij}$, 
such that for all $s\geq s_0$, there exist positive constants 
$\kappa_1(s,\tau)$ and $\kappa$ such that 
$$
\int_{\Omega_T} (s\rho)^{\tau-1}e^{-2s\widetilde{\eta}_{\omega'}}
(\vert \partial_t u\vert^2 + \vert \Delta u\vert^2 
+ (s\rho)^2\vert \nabla u\vert^2 + (s\rho)^4 \vert u\vert^2)
$$
$$
\le \kappa\int_{\Omega_T} (s\rho)^{\tau}
e^{-2s\widetilde{\eta}_{\omega'}}
\vert a_{12}v + A_{12}\cdot\nabla v + f\vert^2
+ \kappa \int_{\omega'_T}(s\rho)^{\tau+3}
e^{-2s\widetilde{\eta}_{\omega'}}\vert u\vert^2
$$
and
$$
\int_{\Omega_T} (s\rho)^{\tau-1}e^{-2s\widetilde{\eta}_{\omega'}}
(\vert \partial_t v\vert^2 + \vert \Delta v\vert^2 
+ (s\rho)^2\vert \nabla v\vert^2 + (s\rho)^4 \vert v\vert^2)
$$
$$
\le \kappa\int_{\Omega_T} (s\rho)^{\tau}
e^{-2s\widetilde{\eta}_{\omega'}}
\vert a_{21}u + A_{21}\cdot\nabla u + g\vert^2
+ \kappa \int_{\omega'_T}(s\rho)^{\tau+3}
e^{-2s\widetilde{\eta}_{\omega'}}\vert v\vert^2
$$
for all large $s > 0$.  Here and henceforth we set
$\widetilde{\eta_{\omega'}}(x,t) = \frac{\beta_{\omega'}(x)}
{t(T-t)}$.  Adding them and choosing $s>0$ sufficiently large
to absorb the terms of $u,v, \nabla u, \nabla v$ on the right hand
side into the left hand side.  
Hence
$$
\int_{\Omega_T} (s\rho)^{\tau-1}e^{-2s\widetilde{\eta}_{\omega'}}
(\vert \partial_t u\vert^2 + \vert \partial_t v\vert^2
+ \vert \Delta u\vert^2 + \vert \Delta v\vert^2
+ (s\rho)^2\vert \nabla u\vert^2 
+ (s\rho)^2\vert \nabla v\vert^2
+ (s\rho)^4 \vert u\vert^2 + (s\rho)^4\vert v\vert^2)
$$
$$
\le \kappa\int_{\Omega_T} (s\rho)^{\tau}e^{-2s\widetilde{\eta}
_{\omega'}}(\vert f\vert^2 + \vert g\vert^2)
+ \kappa\int_{\omega'_T} (s\rho)^{\tau+3}
e^{-2s\widetilde{\eta}_{\omega'}}(\vert u\vert^2 + \vert v\vert^2)
$$
for all large $s>0$.  Moreover we have 
$\vert (s\rho)^{\tau+3}e^{-2s\widetilde{\eta}_{\omega'}}\vert
\le \kappa_2(s,\tau)$ on $\overline{\Omega_T}$ by $\beta_{\omega'}
> 0$ on $\overline{\Omega}$.  Hence
$$
\int_{\omega'_T} (s\rho)^{\tau+3}
e^{-2s\widetilde{\eta}_{\omega'}}(\vert u\vert^2 + \vert v\vert^2)
\le \kappa_2(s,\tau)(\Vert u\Vert^2_{L^2(\omega'_T)}
+ \Vert v\Vert^2_{L^2(\omega'_T)}).
$$
Apply Lemma 2.3, set $\alpha_{\omega} = \beta_{\omega'}$
and note by $\omega' \subset \omega$ that $\Vert u\Vert^2
_{W_2^{2,1}(\omega'_T)} \le \Vert u\Vert^2_{W_2^{2,1}(\omega_T)}$.
Then the proof of Theorem 2.2 is completed.

\section{Proof of Theorem 1.2}
\label{sec: 3}
%
Let us recall that $(U,V)$ satisfie (\ref{syst-UV})  and $(\tilde{U}, \tilde{V})$  system (\ref{syst-UV}) with $b,c,U_0,V_0$  replaced by $\tilde{b},
\tilde{c}, \tilde{U}_0, \tilde{V}_0$ respectively.

We set
\[ 
u = U - \widetilde{U}, \quad v = V - \widetilde{V}. 
\] 
Then $(u,v)$ satisfies 
$$
\partial _{t}u = \Delta u + au + bv + A\cdot\nabla u + B\cdot\nabla v
+ (b-\tilde{b})\tilde{V}, 
$$
$$
\partial _{t}v = \Delta v + cu + dv + C\cdot\nabla u + D\cdot\nabla v
+ (c-\tilde{c})\tilde{U} \quad \mbox{in $\Omega_T$}, 
$$
$$
u = v = 0 \qquad \mbox{on $\Sigma_T$}
$$
and
$$
u(\cdot,\theta) = v(\cdot,\theta) = 0 \qquad \mbox{in $\Omega$}.
$$
By Assumption 1.1, we can assume that 
$\vert \widetilde{U}(x,t)\vert, \vert 
\widetilde{V}(x,t)\vert \ne 0$ for all $(x,t) \in \overline{\Omega_T}$
by taking $T>0$ sufficiently small if necessary.  Moreover
we can assume that $\theta = \frac{T}{2}$. 
Because we take small $\delta > 0$ such that $0\le \theta - \delta
< \theta < \theta + \delta \le T$ and we can replace $\omega
\times (0,T)$ by $\omega \times (\theta-\delta, \theta+\delta)$.
Shifting $t$ by $t-(\theta-\delta)$, we can set $\theta = \delta$
and $T=2\delta$.

Setting 
\[
\widetilde{u} = \frac{u}{\widetilde{V}}, \quad \widetilde{v} = \frac{v}{%
\widetilde{U}},  \quad f = b - \widetilde{b}, \qquad g = c - \widetilde{c},
\]
we have 
\begin{equation}
\partial_t\widetilde{u} = \Delta \widetilde{u} + a_{11}\widetilde{u}
+ a_{12}\widetilde{v}   
+  A_{13}\cdot\nabla \widetilde{u} + A_{14}\cdot \nabla 
\widetilde{v} + f \quad \mbox{in } \Omega _T,       \label{(3.1)}  
\end{equation}
\begin{equation}
\partial_t\widetilde{v} = \Delta \widetilde{v} + a_{21}\widetilde{u}
+ a_{22}\widetilde{v}   
+   A_{23}\cdot\nabla \widetilde{u} + A_{24}\cdot \nabla 
\widetilde{v} + g \quad \mbox{in } \Omega _T,      \label{(3.2)}
\end{equation}
where 
\[
a_{11} = a - \frac{\partial_t\widetilde{V}}{\widetilde{V}} + \frac{%
\Delta \widetilde{V}}{\widetilde{V}} + A\cdot \frac{\nabla \widetilde{V}}{%
\widetilde{V}}, \quad
a_{12} = b\frac{\widetilde{U}}{\widetilde{V}} + B\cdot\frac{\nabla%
\widetilde{U}}{\widetilde{V}}, 
\]
$$
A_{13} = A + \frac{2\nabla \widetilde{V}}{\widetilde{V}},
$$
$$
A_{14}(x,t) = B\frac{\widetilde{U}}{\widetilde{V}}(x,t) \equiv
B(x)W(x,t), 
$$
\[
a_{21} = c\frac{\widetilde{V}}{\widetilde{U}} + C\frac{\nabla 
\widetilde{V}}{\widetilde{U}}, \quad
a_{22} = d - \frac{\partial_t\widetilde{U}}{\widetilde{U}} + \frac{%
\Delta \widetilde{U}}{\widetilde{U}} + D\cdot \frac{\nabla \widetilde{U}}{%
\widetilde{U}} 
\]
and 
\[
A_{23} = C\frac{\widetilde{V}}{\widetilde{U}}, \quad A_{24} = D + 
\frac{2\nabla\widetilde{U}}{\widetilde{U}}. 
\]
Let 
\[
y = \partial_t \widetilde{u}, \quad z = \partial_t\widetilde{v}. 
\]
Since $b,c,\widetilde{b}, \widetilde{c}$ are independent of $t$,
we obtain
$$
\partial_ty = \Delta y + a_{11}y + a_{12}z + A_{13}\cdot\nabla y +
A_{14}\cdot \nabla z   
$$
\begin{equation}
+ (\partial_ta_{11})\widetilde{u} + (\partial_ta_{12})\widetilde{v} +
(\partial_tA_{13})\cdot\nabla \widetilde{u} + (\partial_tA_{14})\cdot 
\nabla \widetilde{v},                \label{(3.3)}
\end{equation}
$$
\partial_tz = \Delta z + a_{21}y + a_{22}z + A_{23}\cdot\nabla y +
A_{24}\cdot \nabla z  
$$
\begin{equation}
+ (\partial_ta_{21})\widetilde{u} + (\partial_ta_{22})\widetilde{v} +
(\partial_tA_{23})\cdot\nabla \widetilde{u} + (\partial_tA_{24})\cdot 
\nabla \widetilde{v}                                \label{(3.4)}
\end{equation}
$$
y = z = 0 \qquad \mbox { on } \Sigma _T.    
$$
{\bf First Step.} In terms of $y$, we estimate an $L^2$-norm of $z$ in 
a subdomain of $\Omega$.  Since $\tilde{u}(x,t) =\int^t_{\theta}
y(x,\xi) d\xi$ and $\tilde{v}(x,t) =\int^t_{\theta} z(x,\xi) d\xi$
by $\tilde{u}(\cdot,\theta) = \tilde{v}(\cdot,\theta) = 0$, we rewrite
  (\ref{(3.3)}) as
$$
B(x)\cdot \nabla z(x,t) + b_1(x)z(x,t) 
+ W_1(x,t)B(x)\cdot\int^t_{\theta} \nabla z(x,\xi)d\xi
+ b_2(x,t)\int^t_{\theta} z(x,\xi) d\xi
$$
$$
=\frac{1}{W(x,t)}\left( \partial_ty(x,t) - \Delta y(x,t) - a_{11}y(x,t)
- A_{13}\cdot\nabla y - (\partial_ta_{11})\int^t_{\theta} y(x,\xi) d\xi
- (\partial_tA_{13})\cdot\int^t_{\theta} \nabla y(x,\xi) d\xi\right)
$$
\begin{equation}
\equiv Q(y)(x,t) \qquad x \in \omega, \thinspace 0<t<T.    \label{(3.5)}
\end{equation}
Here we set 
$$
b_1(x,t) = \frac{a_{12}(x,t)}{W(x,t)}, \quad 
b_2(x,t) = \frac{\partial_ta_{12}(x,t)}{W(x,t)}, \quad
W_1(x,t) = \frac{\partial_t W(x,t)}{W(x,t)}.
$$
We will estimate $z$ in a subdomain $\omega'$ of $\omega$ by means of
 (\ref{(3.5)}), and the argument is similar to Lemma 2.3 but we need a 
special weight function for treating the integral terms 
$\int^t_{\theta} \nabla z(x,\xi)d\xi$ and $\int^t_{\theta} z(x,\xi)
d\xi$.
First we show\\
{\bf Lemma 3.1}  {\it Let $T = 2\theta$ and let $\tilde{  \varphi} \in
C^1[0,T^2]$ and let us assume that there exists a constant $\kappa_0 > 0$
such that $\frac{d\tilde{ \varphi}}{dt}(t) \le - \kappa_0$ for 
$t \in [0, T^2]$.  Then 
$$
\int^T_0 \left\vert \int^t_{\theta} g(\xi)d\xi \right\vert^2
e^{2s\tilde{ \varphi}((t-\theta)^2)} dt
\le \frac{1}{4s\kappa_0}\int^T_0 \vert g(t)\vert^2
e^{2s\tilde{ \varphi}((t-\theta)^2)} dt.
$$
}

The proof is given by Klibanov and Timonov p.78, \cite{KT}.

Henceforth we choose $\tilde{ \varphi}(t) = -t$ and we set 
$ \varphi_1(t) = \tilde{ \varphi}((t-\theta)^2) = -(t-\theta)^2$.  Then the 
conclusion of Lemma 3.1 holds true.

We set 
\begin{equation}
w(x,t) = z(x,t) + W_1(x,t)\int^t_{\theta} z(x,\xi) d\xi, \quad
x\in \Omega, \thinspace 0 < t < T.                \label{(3.6)}
\end{equation}
Then direct calculations yield
\begin{equation}
B(x)\cdot \nabla w(x,t) 
= Q(y)(x,t) - b_1z - b_2\int^t_{\theta} z(x,\xi) d\xi
+ (B\cdot \nabla W_1)\int^t_{\theta} z(x,\xi) d\xi \quad
\mbox{in $\omega_T$}.                                     \label{(3.7)}
\end{equation}
Henceforth $\kappa_j > 0$ denote generic constants which are dependent 
on $M, \delta_0$ in Assumption 1.1 and independent of $s>0$.
In terms of Assumption 1.1 (d), we can apply Lemma 2.3 to obtain
\begin{eqnarray*}
&&s^2\int_{\omega'} \vert w(x,t)\vert^2 e^{2s \varphi_0(x)} dx
\le \kappa_1\int_{\omega'} \vert Q(y)(x,t)\vert^2
e^{2s \varphi_0(x)} dx\\
&+& \kappa_1\int_{\omega'} \vert z(x,t)\vert^2
e^{2s \varphi_0(x)} dx
+ \kappa_1\int_{\omega'} \left\vert \int^t_{\theta} z(x,\xi) d\xi
\right\vert^2 e^{2s \varphi_0(x)} dx
\end{eqnarray*}
for all large $s>0$.  Here and henceforth we set $ \varphi_0(x) = 
x_n - \gamma(x')$.

Hence by Lemma 3.1, we have
\begin{eqnarray*}
&&s^2 \int^T_0\int_{\omega'} \vert w(x,t)\vert^2 
e^{2s( \varphi_0(x)+ \varphi_1(t))} dxdt
\le \kappa_1\int^T_0\int_{\omega'} \vert Q(y)(x,t)\vert^2 
e^{2s( \varphi_0(x)+ \varphi_1(t))} dxdt\\
&+& \kappa_1 \int^T_0\int_{\omega'} \vert z(x,t)\vert^2 
e^{2s( \varphi_0(x)+ \varphi_1(t))} dxdt
+ \kappa_1\int_{\omega'} \left(
\int^T_0 \left\vert \int^t_{\theta} z(x,\xi)d\xi \right\vert^2
e^{2s \varphi_1(t)} dt \right) e^{2s \varphi_0(x)} dx\\
&\le& \kappa_1\int_{\omega'_T} \vert Q(y)(x,t)\vert^2 
e^{2s( \varphi_0(x)+ \varphi_1(t))} dxdt
+ \kappa_1\int_{\omega'_T} \vert z(x,t)\vert^2 
e^{2s( \varphi_0(x)+ \varphi_1(t))} dxdt\\
&+& \frac{\kappa_1}{s}\int_{\omega'_T} \vert z(x,t)\vert^2 
e^{2s( \varphi_0(x)+ \varphi_1(t))} dxdt.
\end{eqnarray*}
Consequently
$$
s^2\int_{\omega'_T} \vert w(x,t)\vert^2 
e^{2s( \varphi_0(x)+ \varphi_1(t))} dxdt
\le \kappa_2e^{2s\kappa_3}\Vert y\Vert^2_{W_2^{2,1}(\omega'_T)}
$$
\begin{equation}
+ \kappa_2\int_{\omega'_T} \vert z(x,t)\vert^2 
e^{2s( \varphi_0(x)+ \varphi_1(t))} dxdt                          \label{(3.8)}
\end{equation}
for all large $s>0$.  On the other hand,  (\ref{(3.6)}) and Lemma 3.1 yield
\begin{eqnarray*}
&&\int_{\omega'_T} \vert z(x,t)\vert^2 
e^{2s( \varphi_0(x)+ \varphi_1(t))} dxdt
= \int_{\omega'_T} \left\vert w(x,t)
- W_1(x,t)\int^t_{\theta} z(x,t) d\xi \right\vert^2 
e^{2s( \varphi_0(x)+ \varphi_1(t))} dxdt\\
&\le& \kappa_4\int_{\omega'_T} \vert w(x,t)\vert^2 
e^{2s( \varphi_0(x)+ \varphi_1(t))} dxdt
+ \kappa_4\int_{\omega'_T} \left\vert 
\int^t_{\theta} z(x,t) d\xi \right\vert^2 
e^{2s( \varphi_0(x)+ \varphi_1(t))} dxdt\\
&\le& \kappa_5 \int_{\omega'_T} \vert w(x,t)\vert^2 
e^{2s( \varphi_0(x)+ \varphi_1(t))} dxdt
+ \frac{\kappa_5}{s}\int_{\omega'_T} 
\vert z(x,t)\vert^2 e^{2s( \varphi_0(x)+ \varphi_1(t))} dxdt
\end{eqnarray*}
for all large $s > 0$.  Hence choosing $s>0$ sufficiently large, we
have
\begin{equation}
\int_{\omega'_T} \vert z(x,t)\vert^2 
e^{2s( \varphi_0(x)+ \varphi_1(t))} dxdt
\le \kappa_6\int_{\omega'_T} \vert w(x,t)\vert^2 
e^{2s( \varphi_0(x)+ \varphi_1(t))} dxdt                        \label{(3.9)}
\end{equation}
for all large $s>0$.  Substituting  (\ref{(3.9)}) into  (\ref{(3.8)}) and fixing
$s>0$ sufficiently large, we obtain
$$
\Vert w\Vert_{L^2(\omega'_T)} \le \kappa_7e^{\kappa_7s}
\Vert y\Vert_{W_2^{2,1}(\omega'_T)}.
$$
Hence by  (\ref{(3.9)}) we have
\begin{equation}
\Vert z\Vert_{L^2(\omega'_T)} \le \kappa_8e^{\kappa_7s}
\Vert y\Vert_{W_2^{2,1}(\omega'_T)}.                    \label{(3.10)}
\end{equation}\\
{\bf Second Step.}
We will estimate $\Vert \nabla z\Vert_{L^2(\omega_1 \times
(\delta,T-\delta))}$ where $\omega_1 \subset \omega$ and $\delta>0$.
For it, we use the interior regularity estimate for a heat equation
 (\ref{(3.4)}) in $z$.  Let us recall that $\rho(t) = \frac{1}{t(T-t)}$.
Setting $\tilde{z}(x,t) = e^{-\rho(t)}z(x,t)$, we rewrite  (\ref{(3.4)}) as
$$
\partial_t\tilde{z}(x,t) = \Delta\tilde{z}(x,t) - \rho'(t)\tilde{z}(x,t)
+ a_{22}\tilde{z} + A_{24}\cdot\nabla \tilde{z}
$$
$$
+ (\partial_ta_{22})\int^t_{\theta} e^{\rho(\xi)-\rho(t)}\tilde{z}(x,\xi)
d\xi 
+ (\partial_tA_{24})\cdot \int^t_{\theta} e^{\rho(\xi)-\rho(t)}
\nabla \tilde{z}(x,\xi)d\xi 
$$
\begin{equation}
+ e^{-\rho(t)}\left(a_{21}y + A_{23}\cdot\nabla y
+ (\partial_ta_{21})\int^t_{\theta} y(x,\xi) d\xi
+ (\partial_tA_{23})\cdot \int^t_{\theta} \nabla y(x,\xi) d\xi\right).
                                        \label{(3.11)}
\end{equation}
We choose subdomains $\omega_1, \omega_2$ of $C^{\infty}$ class
such that $\omega_1 \subset \overline{\omega_1} \subset \omega_2
\subset \overline{\omega_2} \subset \omega'$ and choose $\chi 
\in C^1(\overline{\omega'})$, $\ge 0$ such that 
$$
\chi(x)=
\left\{
\begin{array}{rl} 
1, &\quad x \in \omega_1,\\
0, &\quad x \in \omega' \setminus \overline{\omega_2}.
\end{array}\right.
$$
Moreover we can take $\chi$ satisfying
\begin{equation}
\frac{\vert \nabla\chi(x)\vert^2}{\chi(x)} \le \kappa_9 
\quad x\in \overline{\omega'}                            \label{(3.12)}
\end{equation}
(e.g., p.414 in Lions \cite{Li}).
Multiplying  (\ref{(3.11)}) with $\chi\tilde{z}$ and integrating over 
$\omega' \times (0,T)$, we have
\begin{eqnarray*}
&& \frac{1}{2}\int^T_0\int_{\omega'}
\chi(x)\partial_t(\tilde{z}^2) dxdt
= -\int^T_0 \int_{\omega'} \chi\vert \nabla \tilde{z}\vert^2 dxdt
- \int^T_0 \int_{\omega'} \nabla\chi \cdot \tilde{z}\nabla\tilde{z} dxdt
- \int^T_0 \int_{\omega'} \chi \rho'(t)e^{-2\rho(t)}\vert z\vert^2 dxdt\\
&+& \int^T_0 \int_{\omega'} (a_{22}\vert \tilde{z}\vert^2\chi
+ A_{24}\cdot\nabla \tilde{z}\chi\tilde{z}) dxdt
+   \int^T_0 \int_{\omega'}
(\partial_t a_{22})\chi\tilde{z} \left( \int^t_{\theta} e^{\rho(\xi)-\rho(t)}
\tilde{z}(x,\xi) d\xi\right) dxdt\\
&+&   \int^T_0 \int_{\omega'}(\partial_t A_{24}) \cdot
\chi\tilde{z} \left( \int^t_{\theta} e^{\rho(\xi)-\rho(t)}
\nabla \tilde{z}(x,\xi) d\xi\right) dxdt\\
&+& \int^T_0\int_{\omega'} e^{-\rho(t)}\chi\tilde{z}
\left( a_{21}y + A_{23}\nabla y
+ (\partial_ta_{21})\int^t_{\theta} y(x,\xi) d\xi
+ (\partial_tA_{23})\cdot \int^t_{\theta} \nabla y(x,\xi) d\xi \right)
dxdt.
\end{eqnarray*}
By the Cauchy-Schwarz inequality and  (\ref{(3.12)}), we have
$$
\vert \nabla \chi\cdot \tilde{z}\nabla\tilde{z}\vert
= \left\vert \frac{\nabla\chi}{\sqrt{\chi}}\tilde{z}\cdot
\sqrt{\chi}\nabla \tilde{z} \right\vert
\le \frac{1}{8}\chi\vert \nabla \tilde{z}\vert^2 
+ \frac{2\vert \nabla\chi\vert^2}{\chi}\vert \tilde{z}\vert^2
$$
and
$$
\vert A_{24}\cdot \chi\tilde{z}\nabla\tilde{z}\vert
= \vert \sqrt{\chi}\nabla\tilde{z}\cdot A_{24}\sqrt{\chi}\tilde{z} \vert
\le \frac{1}{8}\chi\vert \nabla \tilde{z}\vert^2 
+ 2\vert A_{24}\vert^2\chi\vert \tilde{z}\vert^2.
$$
Hence, since $\tilde{z}(\cdot,0) = \tilde{z}(\cdot,T) = 0$, 
$$
\sup_{0\le t\le T} \vert \rho'(t)e^{-2\rho(t)} \vert < \infty
$$
and $\rho(\xi) - \rho(t) \le 0$ if $\xi$ is between $\theta$ and 
$t$, we have
\begin{eqnarray*}
&&\int^T_0\int_{\omega'} \chi\vert \nabla\tilde{z}\vert^2 dxdt
\le \frac{1}{4}\int^T_0\int_{\omega'} \chi\vert \nabla\tilde{z}\vert^2 
dxdt\\
&+& \kappa_{10}\int^T_0\int_{\omega'} (\vert \tilde{z}\vert^2
+ \vert z\vert^2) dxdt
+ \kappa_{10}\int^T_0\int_{\omega'} \vert \tilde{z}\vert
\left\vert \int^t_{\theta} \tilde{z}(x,\xi) d\xi\right\vert dxdt\\
&+& \kappa_{10}\int^T_0\int_{\omega'} \chi\vert \tilde{z}\vert
\left\vert \int^t_{\theta} \nabla \tilde{z}(x,\xi) d\xi\right\vert dxdt\\
&+& \kappa_{10}\int^T_0\int_{\omega'} 
\left\{ \vert \tilde{z}\vert(\vert y\vert + \vert \nabla y\vert)
+ \vert \tilde{z}\vert
\left( \left\vert \int^t_{\theta} y(x,\xi) d\xi\right\vert 
+ \left\vert \int^t_{\theta} \nabla y(x,\xi) d\xi\right\vert \right)
\right\} dxdt.
\end{eqnarray*}
Moreover the Cauchy-Schwarz inequality yields 
\begin{eqnarray*}
&& \int^T_0\int_{\omega'} \chi\vert \tilde{z}(x,t)\vert
\left\vert \int^t_{\theta} \nabla\tilde{z}(x,\xi) d\xi \right\vert dxdt\\
&\le& \int^T_0\int_{\omega'} \sqrt{\chi}\vert \tilde{z}(x,t)\vert
\left( \int^T_0 \sqrt{\chi(x)} \vert \nabla \tilde{z}(x,\xi) \vert d\xi
\right) dxdt\\
&\le& \int^T_0\int_{\omega'} \left( \frac{1}{8T^2}
\left\vert \int^T_0 \sqrt{\chi}\vert \nabla \tilde{z}(x,\xi)\vert
d\xi \right\vert^2 + 2T^2\chi \vert \tilde{z}(x,\xi) \vert^2\right)
dxdt\\
&\le& 2T^2\int^T_0\int_{\omega'} \vert \tilde{z}(x,t)\vert^2 dxdt
+ \int^T_0\int_{\omega'} \frac{1}{8T}\int^T_0
\chi \vert \nabla \tilde{z}(x,\xi) \vert^2d\xi dxdt\\
&\le& 2T^2\int^T_0\int_{\omega'} \vert \tilde{z}(x,t)\vert^2 dxdt
+ \frac{1}{8}\int^T_0\int_{\omega'} \chi 
\vert \nabla \tilde{z}(x,t) \vert^2 dxdt.
\end{eqnarray*}
Hence
\begin{eqnarray*}
&&\frac{5}{8}\int^T_0\int_{\omega'} \chi 
\vert \nabla \tilde{z}(x,t) \vert^2 dxdt\\
&\le& \kappa_{11}\int^T_0\int_{\omega'}
(\vert \tilde{z}\vert^2 + \vert z\vert^2) dxdt
+ \kappa_{11}\int^T_0\int_{\omega'}
(\vert y\vert^2 + \vert \nabla y\vert^2) dxdt.
\end{eqnarray*}
Let $\delta > 0$ be fixed sufficiently small.  Then 
$\vert \nabla \tilde{z}(x,t)\vert \ge \kappa_{12}(\delta)
\vert \nabla z(x,t)\vert$ for $\delta \le t \le T-\delta$.
Since $\chi=1$ in $\omega_1$, we have
$$
\int^{T-\delta}_{\delta} \int_{\omega_1}
\vert \nabla z\vert^2 dxdt
\le \kappa_{13}(\delta)(\Vert z\Vert^2_{L^2(\omega_T')}
+ \Vert y\Vert^2_{L^2(0,T;H^1(\omega'))}).
$$
By means of  (\ref{(3.10)}), we obtain
\begin{equation}
\Vert z\Vert_{L^2(\delta,T-\delta;H^1(\omega_1))}
\le \kappa_{14}(\delta)\Vert y\Vert_{W_2^{2,1}(\omega'_T)}.
                                                          \label{(3.13)}
\end{equation}
{\bf Third Step.}
We apply Theorem 2.2 to   (\ref{(3.3)}) and  (\ref{(3.4)}) for $\omega' \subset 
\Omega$ and $(\delta, T-\delta)$.  We set
$$
\eta(x,t) = \frac{\alpha_{\omega'}(x)}{(t-\delta)(T-\delta-t)}.
$$
Using also  (\ref{(3.13)}), we obtain that there exist two positive constants 
$s_{0}$ and $\kappa $  such that for all $s\geq s_0$, one has
$$
\int^{T-\delta}_{\delta} \int_{\Omega} (s\rho)^{-1}e^{-2s\eta}
(\vert \partial_t y\vert^2 + \vert \partial_t z\vert^2
+ \vert \Delta y\vert^2 + \vert \Delta z\vert^2
+ (s\rho)^2\vert \nabla y\vert^2 
+ (s\rho)^2\vert \nabla z\vert^2
+ (s\rho)^4 \vert y\vert^2 + (s\rho)^4\vert z\vert^2)dxdt 
$$
$$
\le \kappa_{15}\int^{T-\delta}_{\delta} \int_{\Omega}
\vert (\partial_ta_{11})\widetilde{u} 
+ (\partial_ta_{12})\widetilde{v} 
+ (\partial_tA_{13})\cdot\nabla\widetilde{u} 
+ (\partial_tA_{14})\cdot\nabla\widetilde{v}\vert^2  e^{-2s\eta}dxdt
$$
$$
+ \kappa_{15} \int^{T-\delta}_{\delta} \int_{\Omega}
\vert (\partial_ta_{21})\widetilde{u} 
+ (\partial_ta_{22})\widetilde{v} 
+ (\partial_tA_{23})\cdot \nabla\widetilde{u} 
+ (\partial_tA_{24})\cdot \nabla\widetilde{v}\vert^2)e^{-2s\eta}dxdt
$$
$$
+ \kappa_{16}(s)\Vert y\Vert^2_{W_2^{2,1}(\omega'\times (\delta,
T-\delta))}
$$
$$
+ \kappa_{16}(s)\Biggl\Vert
(\partial_ta_{11})\int^t_{\theta} y(x,\xi) d\xi
+ (\partial_ta_{12})\int^t_{\theta} z(x,\xi) d\xi  
+ (\partial_tA_{13})\cdot \int^t_{\theta} \nabla y(x,\xi) d\xi
$$
$$
+ (\partial_tA_{14})\cdot \int^t_{\theta} \nabla z(x,\xi) d\xi
\Biggr\Vert^2_{L^2(\omega' \times (\delta,T-\delta))}
$$
$$
+ \kappa_{16}(s)\Biggl\Vert
(\partial_ta_{21})\int^t_{\theta} y(x,\xi) d\xi
+ (\partial_ta_{22})\int^t_{\theta} z(x,\xi) d\xi  
+ (\partial_tA_{23})\cdot \int^t_{\theta} \nabla y(x,\xi) d\xi
$$
$$
+ (\partial_tA_{24})\cdot \int^t_{\theta} \nabla z(x,\xi) d\xi
\Biggr\Vert^2_{L^2(\omega' \times (\delta,T-\delta))}
$$
$$
\le \kappa_{15}\int^{T-\delta}_{\delta} \int_{\Omega}
\vert (\partial_ta_{11})\widetilde{u} 
+ (\partial_ta_{12})\widetilde{v} 
+ (\partial_tA_{13})\cdot\nabla\widetilde{u} 
+ (\partial_tA_{14})\cdot\nabla\widetilde{v}\vert^2  e^{-2s\eta}dxdt
$$
$$
+ \kappa_{15} \int^{T-\delta}_{\delta} \int_{\Omega}
\vert (\partial_ta_{21})\widetilde{u} 
+ (\partial_ta_{22})\widetilde{v} 
+ (\partial_tA_{23})\cdot \nabla\widetilde{u} 
+ (\partial_tA_{24})\cdot \nabla\widetilde{v}\vert^2)e^{-2s\eta}dxdt
$$
\begin{equation}
+ \kappa_{16}(s)\Vert y\Vert^2_{W_2^{2,1}(\omega_T)} .   \label{(3.14)}
\end{equation}
for all large $s > 0$.
In order to improve inequality  (\ref{(3.14)}), we use the following 
lemma. (\cite{KT:04}~, Lemma 3.1.1 in \cite{KT}).
\\
{\bf Lemma 3.2}
{\it Let $\theta = \frac{T}{2}$. There exists a positive constant 
$\kappa_{17}$ such that 
\[
\int^{T-\delta}_{\delta} \int_{\Omega}
\left\vert \int^t_{\theta} q(x,\xi) d\xi\right\vert^2 
e^{-2s\eta} dxdt
\le \frac{\kappa_{17}}{s}\int^{T-\delta}_{\delta}\int_{\Omega} 
\vert q(x,t)\vert^2 e^{-2s\eta}dxdt
\]
for $s > 0$.}
\\
{\bf Proof of Lemma 3.2} The proof is similar to \cite{KT:04}, 
Lemma 3.1.1 in \cite{KT}.  We have
$$
\int^{T-\delta}_{\delta} \int_{\Omega}
\left\vert \int^t_{\theta} q(x,\xi) d\xi\right\vert^2 
e^{-2s\eta} dxdt
$$
$$
= \int_{\Omega} \int^{\theta}_{\delta} 
\left\vert \int^{\theta}_t q(x,\xi) d\xi\right\vert^2 
e^{-2s\eta} dxdt
+ \int_{\Omega} \int^{T-\delta}_{\theta} 
\left\vert \int^t_{\theta} q(x,\xi) d\xi\right\vert^2 
e^{-2s\eta} dxdt.
$$
It is sufficient to estimate the second term because the estimation of
the first term is similar.
By the Cauchy-Schwarz inequality, we obtain
$$
\int_{\Omega} \int^{T-\delta}_{\theta} 
\left\vert \int^t_{\theta} q(x,\xi) d\xi\right\vert^2 
e^{-2s\eta} dxdt
\le \int_{\Omega} \int^{T-\delta}_{\theta} 
(t-\theta)\left( \int^t_{\theta} \vert q(x,\xi)\vert^2 d\xi\right) 
e^{-2s\eta} dxdt
$$
$$
\le \int_{\Omega} \int^{T-\delta}_{\theta} 
\frac{(t-\delta)^2(T-\delta-t)^2}{2\alpha_{\omega'}(x)}
\partial_t\eta(x,t)\left(\int^t_{\theta} \vert q(x,\xi)\vert^2 
d\xi\right)e^{-2s\eta} dxdt.
$$
Here we used
$$
\partial_t\eta(x,t)
= \frac{2(t-\theta)\alpha_{\omega'}(x)}
{(t-\delta)^2(T-\delta-t)^2}.
$$
Noting that $\alpha_{\omega'}(x) > 0$ and $\partial_t\eta(x,t)
\ge 0$ for $x \in \overline{\Omega}$ and $\theta \le t \le T-\delta$,
we have
$$
\int_{\Omega} \int^{T-\delta}_{\theta} 
\left( \int^t_{\theta} \vert q(x,\xi)\vert^2 d\xi\right)
e^{-2s\eta}dxdt
\le \kappa'_{17}
\int_{\Omega} \int^{T-\delta}_{\theta} 
\left( \int^t_{\theta} \vert q(x,\xi)\vert^2 d\xi\right)
(\partial_t\eta(x,t))e^{-2s\eta(x,t)} dxdt
$$
$$
= -\frac{\kappa'_{17}}{2s}
\int_{\Omega} \int^{T-\delta}_{\theta} 
\left( \int^t_{\theta} \vert q(x,\xi)\vert^2 d\xi\right)
\partial_t(e^{-2s\eta(x,t)}) dxdt.
$$
By noting that $e^{-2s\eta(x,T-\delta)}=0$, the integration by parts
implies that the right hand side is equal to 
$$
\frac{\kappa'_{17}}{2s}
\int_{\Omega} \int^{T-\delta}_{\theta} 
\vert q(x,t)\vert^2 \partial_t(e^{-2s\eta(x,t)}) dxdt.
$$
Thus the proof of Lemma 3.2 is completed.

Since $\widetilde{u}(x,t) = \int^t_{\theta} y(x,\xi)d\xi$ and
$\widetilde{v}(x,t) = \int^t_{\theta} z(x,\xi)d\xi$,
by a direct application of this lemma,  the first integral on the 
right hand side of  (\ref{(3.14)}) can be absorbed into the left hand 
side.
Hence 
\begin{eqnarray*}
&&\int^{T-\delta}_{\delta}\int_{\Omega} 
\Biggl\{ \frac{1}{s\rho}(\vert \partial_ty\vert^2 
+ \vert \partial_tz\vert^2 + \vert \Delta y\vert^2 
+ \vert \Delta z\vert^2) 
+ s\rho(\vert \nabla y\vert^2 + \vert \nabla z\vert^2) 
+ s^3\rho^3(\vert y\vert^2 + \vert z\vert^2)\Biggr\} e^{-2s\eta}dxdt  \\
&\leq& \kappa_{16}(s)
(\Vert \partial_t(U - \widetilde{U})\Vert^2_{W_2^{2,1}(\omega_T)}
+ \Vert U - \widetilde{U}\Vert^2_{W_2^{2,1}(\omega_T)})
\end{eqnarray*}
for all large $s > 0$.
We choose $t_0 > 0$ sufficiently small such that $\delta <
t_0  < \theta < T-t_0 < T-\delta$, so that 
\begin{eqnarray*}
&&\int^{T-t_0}_{t_0} \int_{\Omega} \Biggl\{ \frac{1}{s\rho}(\vert
\partial_ty\vert^2 + \vert \partial_tz\vert^2 + \vert \Delta y\vert^2 +
\vert \Delta z\vert^2) + s\rho(\vert \nabla y\vert^2 + \vert \nabla
z\vert^2) + s^3\rho^3(\vert y\vert^2 + \vert z\vert^2)\Biggr\}
e^{-2s\eta} dxdt  \\
&\leq& \kappa_{16}(s) 
(\Vert \partial_t(U - \widetilde{U})\Vert^2_{W_2^{2,1}(\omega_T)}
+ \Vert U - \widetilde{U}\Vert^2_{W_2^{2,1}(\omega_T)}).
\end{eqnarray*}
Since $\frac{1}{\rho}e^{-2s\eta}$, $\rho e^{-2s\eta} \ge \kappa_0(t_0,s)$
on $\overline{\Omega} \times [t_0, T-t_0]$, we fix $s > 0$ sufficiently
large, so that 
\begin{eqnarray*}
&&\Vert \widetilde{u}\Vert^2_{H^1(t_0,T-t_0;H^2(\Omega))} 
+ \Vert \widetilde{u}\Vert^2_{H^2(t_0,T-t_0;L^2(\Omega))}
+ \Vert \widetilde{v}\Vert^2_{H^1(t_0,T-t_0;H^2(\Omega))} 
+ \Vert \widetilde{v}\Vert^2_{H^2(t_0,T-t_0;L^2(\Omega))}\\
&\leq& \kappa_{16}(s)
(\Vert \partial_t(U - \widetilde{U})\Vert^2_{W_2^{2,1}(\omega_T)}
+ \Vert U - \widetilde{U}\Vert^2_{W_2^{2,1}(\omega_T)}).
\end{eqnarray*} 
By the trace theorem, we have 
\begin{eqnarray*}
&&\Vert \partial_t\widetilde{u}(\cdot,\theta)\Vert^2_{L^2(\Omega)} + \Vert
\partial_t\widetilde{v}(\cdot,\theta)\Vert^2_{L^2(\Omega)}
+ \Vert \widetilde{u}(\cdot,\theta)\Vert^2_{H^2(\Omega)}
+ \Vert \widetilde{v}(\cdot,\theta)\Vert^2_{H^2(\Omega)}\\
&\leq& \kappa_{16}(s)
(\Vert \partial_t(U - \widetilde{U})\Vert^2_{W_2^{2,1}(\omega_T)}
+  \Vert U - \widetilde{U}\Vert^2_{W_2^{2,1}(\omega_T)}).
\end{eqnarray*}
Since $f$ and $g$ satisfy  (\ref{(3.1)}) and  (\ref{(3.2)}) at $t=\theta$, we see that
$$
\Vert b-\widetilde{b}\Vert^2_{L^2(\Omega)} + \Vert c-\widetilde{c}
\Vert^2_{L^2(\Omega)} 
\leq \kappa 
(\Vert \partial_t(U - \widetilde{U})\Vert^2_{W_2^{2,1}(\omega_T)}
+ \Vert U - \widetilde{U}\Vert^2_{W_2^{2,1}(\omega_T)}).
$$
Thus the proof of Theorem \ref{theo06} is completed. 
%
\section{Removing the positivity assumption}\label{sec: 4}
%
\noindent
For the stability in our inverse problem, the non-vanishing condition 
Assumption 1.1 (f) is crucial and does not hold automatically. 
We are going to prove that one can realize this assumption 
by a suitable control. 

Let $m \in {\bf N}$ be fixed such that
\begin{equation}
\frac{m}{4} > n.                          \label{(4.1)}
\end{equation}
We assume that 
\begin{equation}
a, \widetilde{b}, \widetilde{c}, d, A, B, C, D 
\in W^{2m-2,\infty}(\Omega).                       \label{(4.2)}
\end{equation}
We set 
$$
L(u,v) = L(a,\widetilde{b},\widetilde{c},d,A,B,C,D)(u,v) =\left( 
\begin{array}{cc}
L_1(u,v) \\ 
L_2(u,v) \\
\end{array}
\right)
$$
\begin{equation}
= -\left( 
\begin{array}{cc}
\Delta u + au + \widetilde{b}v + A\cdot \nabla u +  B\cdot \nabla v\\ 
\Delta v + \widetilde{c}u + dv + C\cdot \nabla v + D\cdot \nabla v \\
\end{array}
\right)                               \label{(4.3)}
\end{equation}
and
$$
D(L)=\left( H^{2}(\Omega )\cap H_{0}^{1}(\Omega )\right) ^{2}.
$$
For $h \in L^2(\omega_T)$, let $(\widetilde{U},\widetilde{V})
:=(\widetilde{U}(\tilde{U}_0,\tilde{V}_0, h)(\cdot,\cdot), 
\widetilde{V}(\tilde{U}_0,\tilde{V}_0,h)(\cdot,\cdot))$
satisfy 
$$
\partial_{t}(\widetilde{U},\widetilde{V}) 
= -L(a,\widetilde{b},\widetilde{c},d,A,B,C,D)
(\widetilde{U}, \widetilde{V}) + (\chi_{\omega}h,0) \quad
\mbox{ \ in }\Omega _T,   
$$
\begin{equation}
(\widetilde{U}, \widetilde{V}) = (0,0) \thinspace\mbox{ on }
\Sigma_T, \quad (\tilde{U},\tilde{V})(\cdot,0) = (\tilde{U}_0,\tilde{V}_0)
\thinspace \mbox{in $\Omega$}.                        \label{(4.4)}
\end{equation}
By $(U,V)$ we denote the solution to  (\ref{(4.4)}) with $b,c$ replacing 
$\widetilde{b},\widetilde{c}$.  Our main result in this section is 
the following~:
\\
{\bf Theorem 4.1}
{\it Suppose Assumption 1.1 except for (f). 
Let $\omega_1$ be a neighbourhood of $\partial \Omega$ such that 
$\overline{\omega} \subset \omega_1$ and let 
$b = \widetilde{b}$ and $c = \widetilde{c}$ in $\omega_1$. 
Let $(U,V)(\cdot,\theta )=(\widetilde{U},\widetilde{V})(\cdot,\theta)$.
Then there exists $h\in L^2(\omega_T)$ depending on 
$a,\widetilde{b},\widetilde{c},
d,A,B,C,D$, $\widetilde{U_0},\widetilde{V_0}$ and $\omega$, 
such that there exists a constant $\kappa > 0$ such that
\begin{equation}
\Vert b - \widetilde{b}\Vert_{L^2(\Omega)}
+ \Vert c - \widetilde{c}\Vert_{L^2(\Omega)}
\le \kappa 
(\Vert \partial_t(U - \widetilde{U})\Vert_{ W_2^{2,1}(\omega_T)}
+ \Vert U - \widetilde{U}\Vert_{ W_2^{2,1}(\omega_T)})      \label{(4.5)}
\end{equation}
for arbitrary $b,c,U,V$ satisfying Assumption 1.1 (a), (e), (h).
}

The rest of this section is devoted to the proof of Theorem 4.1.

{\bf First Step.}
First we prove
\\
{\bf Lemma 4.2} 
{\it 
Let Assumption 1.1 except for (f) hold and let $b = \tilde{b}$,
$c = \tilde{c}$ in $\omega_1$.  Then there exists $h \in L^2(\omega_T)$
such that
\begin{equation}
\vert \tilde{U}(\cdot,\theta)\vert, \thinspace
\vert \tilde{V}(\cdot,\theta)\vert \ne 0 \quad \mbox{on 
$\overline{\Omega \setminus \omega_1}$}.                  \label{(4.6)}
\end{equation}
}

In this step, we will give the proof of Lemma 4.2, which is 
based on the approximate controllability and our Carleman estimate
Theorem 2.2.

Taking $M>0$ for $a,\widetilde{b},
\widetilde{c}, d, A,B,C,D$, and setting $U_1 = e^{-Mt}\widetilde{U}$
and $V_1 = e^{-Mt}\widetilde{V}$, we have
$$
\partial_tU_1 = \Delta U_1 + (a-M)U_1 + \widetilde{b}V_1
+ A\cdot\nabla U_1 + B\cdot \nabla V_1 + e^{-Mt}\chi_{\omega}h
$$
and
$$
\partial_tV_1 = \Delta V_1 + \widetilde{c}U_1 + (d-M)V_1
+ C\cdot\nabla U_1 + D\cdot \nabla V_1.
$$
Consequently, by choosing $M>0$ sufficiently large, 
the integration by parts yields
$$
((L+MI)(u,v), (u,v))_{(L^2(\Omega))^2}
\ge \kappa_1\Vert (u,v)\Vert^2_{(H^1(\Omega))^2}, \quad
(u,v) \in {\mathcal{D}}(L).
$$
Therefore with fear of confusion, we may denote $a-M$ and $d-M$ by 
$a$ and $d$ respectively.  Then 
\begin{equation}
\Vert (u,v)\Vert_{(H^1(\Omega))^2} 
\le \kappa_1\Vert L(u,v)\Vert_{(L^2(\Omega))^2}, \quad
(u,v) \in {\mathcal{D}}(L).                               \label{(4.7)}
\end{equation}
Here and henceforth $\kappa_j>0$ denote generic constants which 
depend on $\Omega$, $M$, $\Vert a\Vert_{W^{2m-2,\infty}(\Omega)}$,
$\Vert \widetilde{b}\Vert_{W^{2m-2,\infty}(\Omega)}$,
$\Vert \widetilde{c}\Vert_{W^{2m-2,\infty}(\Omega)}$,
$\Vert d\Vert_{W^{2m-2,\infty}(\Omega)}$,
$\Vert A\Vert_{(W^{2m-2,\infty}(\Omega))^n}$,
$\Vert B\Vert_{(W^{2m-2,\infty}(\Omega))^n}$,
$\Vert C\Vert_{(W^{2m-2,\infty}(\Omega))^n}$,
$\Vert D\Vert_{(W^{2m-2,\infty}(\Omega))^n}$.
We can prove
\\
{\bf Lemma 4.3}
{\it  
Under assumption  (\ref{(4.2)}), there exists a constant $\kappa_2 > 0$
such that
$$
\Vert (u,v)\Vert_{(H^{2m}(\Omega))^2} 
\le \kappa_2\Vert L^m(u,v)\Vert_{(L^2(\Omega))^2}, \quad
(u,v) \in {\mathcal{D}}(L^m).                               
$$
}
\\
{\bf Proof of Lemma 4.3}
The proof is done by the classical regularity property 
for the Dirichlet problem for the Poisson equation
(e.g., Theorem 8.13 in Gilbarg and Trudinger \cite{GT}) and 
given here for completeness.

We recall  (\ref{(4.3)}) and we set $Q(u,v) = (Q_1(u,v), Q_2(u,v))$,
$Q_1(u,v) = au + \tilde{b}v + A\cdot\nabla u + B\cdot\nabla v$
and $Q_2(u,v) = \tilde{c}u + dv + C\cdot\nabla u + D\cdot\nabla v$.
Let $(u,v) \in {\mathcal{D}}(L^m)$.  By the elliptic regularity 
(e.g., Theorem 8.13 in \cite{GT}) in the 
Dirichlet problem for $\Delta u = f$, we have
$$
\Vert u\Vert_{H^2(\Omega)} \le \kappa_1
\Vert (-L_1-Q_1)(u,v)\Vert_{L^2(\Omega)},
$$
and
$$
\Vert v\Vert_{H^2(\Omega)} \le \kappa_1
\Vert (-L_2-Q_2)(u,v)\Vert_{L^2(\Omega)},
$$
so that
$$
\Vert (u,v)\Vert_{(H^2(\Omega))^2} 
\le \kappa_2\Vert L(u,v)\Vert_{(L^2(\Omega))^2}
+ \kappa_2\Vert (u,v)\Vert_{(H^1(\Omega))^2}.
$$
Hence by  (\ref{(4.7)}), we have
\begin{equation}
\Vert (u,v)\Vert_{(H^2(\Omega))^2} 
\le \kappa_3\Vert L(u,v)\Vert_{(L^2(\Omega))^2}.    \label{(4.8)}
\end{equation}
Again the elliptic regularity yields
$$
\Vert (u,v)\Vert_{(H^3(\Omega))^2} 
\le \kappa_1\Vert (\Delta u, \Delta v)\Vert_{(H^1(\Omega))^2}
+ \kappa_1\Vert (u,v)\Vert_{(L^2(\Omega))^2} 
$$
$$
\le \kappa_1\Vert (-L-Q)(u,v)\Vert_{(H^1(\Omega))^2}
+ \kappa_1\Vert (u,v)\Vert_{(L^2(\Omega))^2} 
$$
\begin{equation}
\le \kappa_1 \Vert L(u,v)\Vert_{(H^1(\Omega))^2}
+ \kappa_1\Vert (u,v) \Vert_{(H^2(\Omega))^2}.             \label{(4.9)}
\end{equation}
On the other hand, we have $L(u,v) \in {\mathcal{D}}(L)$ and apply
 (\ref{(4.7)}) to $L(u,v)$ to have 
$$
\Vert L(u,v)\Vert_{(H^1(\Omega))^2} 
\le \kappa_1\Vert L^2(u,v)\Vert_{(L^2(\Omega))^2}.
$$
Applying this and  (\ref{(4.8)}) to  (\ref{(4.9)}), we obtain
$$
\Vert (u,v)\Vert_{(H^3(\Omega))^2} 
\le \kappa_4\Vert L^2(u,v)\Vert_{(L^2(\Omega))^2}.
$$
Repeating these arguments, we can complete the proof of 
Lemma 4.3.

Moreover by \cite{Pa} and \cite {T} for example, we see:
\\
{\bf Lemma 4.4}
{\it The operator $-L$ generates an analytic semigroup in 
$(L^2(\Omega))^2$. 
}
\\

There are no general result on the approximate controllabilty
for parabolic systems with controls of a restricted number of
components and see e.g., \cite {ABD} and \cite{Leiva} as
related works.  For controllability for systems, see 
\cite {ABD} - \cite{ABDK}, \cite{GPe1} - \cite{GdT}.
Next we will prove the approximate controllability with control
$\chi_{\omega}h$ to only one component.

\noindent
{\bf Lemma 4.5}
{\it
For any $ \varepsilon >0$, $(\tilde{U}_0, \tilde{V}_{0})\in 
(L^{2}(\Omega))^2$, $(\tilde{U}_{1},\tilde{V}_{1})\in (L^{2}(\Omega))^{2}$,
and any $t_0 \in (0, \delta)$, there exists $h_{ \varepsilon}\in 
L^{2}(\omega_T)$ such that
$$
\Vert \tilde{U}(\tilde{U}_0,\tilde{V}_0,h)(\cdot,t_0)
- \tilde{U}_1\Vert_{L^2(\Omega)}
+ \Vert \tilde{V}(\tilde{U}_0,\tilde{V}_0,h)(\cdot,t_0)
- \tilde{V}_1\Vert_{L^2(\Omega)} <  \varepsilon.
$$
}
\\
{\bf Proof.} 
Consider the following reaction-diffusion-convection system : 
$$
\partial _{t}u=\Delta u+au+ \tilde{b}v-\nabla \cdot(Au)
- \nabla \cdot (Bv) \quad \mbox{in $\Omega _{T}$},
$$
$$ 
\partial _{t}v=\Delta v+ \tilde{c}u + dv - \nabla \cdot (Cu)
- \nabla \cdot (Dv) \quad \mbox{in $\Omega _{T}$}, 
$$
\begin{equation}
u=v=0 \qquad \mbox{on $\Sigma _{T}$}.        \label{(4.10)}
\end{equation}
The approximate controllability is equivalent to the uniqueness:
Let $u, v$ satisfy  (\ref{(4.10)}).  Then $u=0$ in $\omega_T$ implies
$u=v=0$ in $\Omega_T$ (e.g., Zabczyk \cite{Z}).  This uniqueness 
follows from Theorem 2.2 by replacing the coefficients in (\ref{syst-uv}) 
suitably and verifying Assumption 1.1 (d).

Now we will complete
\\
{\bf Proof of Lemma 4.2}
The proof is be done in three steps.
Henceforth for fixed $(\widetilde{U}_0, \widetilde{V}_0)$, by 
$(\widetilde{U}, \widetilde{V})(h)$ we denote
$(\widetilde{U}, \widetilde{V})(\tilde{U}_0, \tilde{V}_0,h)$.
\\
{\bf Existence of a control in $L^2(\omega _T)$}\\
Let us arbitrarily fix $(\tilde{U}_1, \tilde{V}_1) \in 
(H^{2m+2}_0(\Omega))^2$ satisfying $\vert \tilde{U}_1\vert$, 
$\vert \tilde{V}_1\vert \ne 0$ on 
$\overline{\Omega \setminus \omega_1}$. 
Then for any $ \varepsilon > 0$ and any $T_1 \in (0, \theta)$,
there exists $h_{ \varepsilon} \in L^2(\omega_{T_1})$ such that
\begin{equation}
\Vert(\widetilde{U}, \widetilde{V})(h_{ \varepsilon})(\cdot,T_1) 
- (\tilde{U}_1, \tilde{V}_1)
\Vert_{(L^2(\Omega))^2}\leq  \varepsilon.         \label{(4.11)}
\end{equation}
\\
{\bf A more regular control}\\
By the density of $C^{\infty}(\omega_{T_1})$ in $L^2(\omega_{T_1})$, 
for any $\delta >0$, there exists $h_{ \varepsilon,\delta}
\in C^{\infty}(\omega_{T_1})$ such that
\begin{equation}
\Vert h_{ \varepsilon}-h_{ \varepsilon,\delta}\Vert
_{L^2(\omega _{T_1})}\leq \delta.                      \label{(4.12)}
\end{equation}
Therefore
$$
\Vert (\widetilde{U},\widetilde{V})(h_{ \varepsilon})(\cdot,T_1)
- (\widetilde{U},\widetilde{V})(h_{ \varepsilon,\delta})
(\cdot,T_1)\Vert_{(L^2(\Omega))^2} 
= \left\Vert \int^{T_1}_0 e^{-(T_1-s)L}\chi_{\omega}
(h_{ \varepsilon} - h_{ \varepsilon,\delta})(s) ds \right\Vert
_{(L^2(\Omega))^2} \le \kappa_5\delta.
$$
\\
{\bf Use of the time regularizing effect}
\\
By  (\ref{(4.11)}) and  (\ref{(4.12)}), we obtain
\begin{equation}
\Vert (\widetilde{U},\widetilde{V})(h_{ \varepsilon,\delta})
(\cdot,T_1) - (\widetilde{U}_1, \widetilde{V}_1)\Vert_{(L^2(\Omega))^2}
\le  \varepsilon + \kappa_5\delta.                   \label{(4.13)}           \end{equation}
Since $-L$ generates an analytic semigroup in $(L^2(\Omega))^2$, 
by e.g., \cite {Pa}, \cite{T}, we see that 
$$
e^{-(\theta-T_1)L}
(\widetilde{U},\widetilde{V})(h_{ \varepsilon,\delta})(\cdot,T_1)
\in D(L^m)
$$
and
$$
\Vert L^m[e^{-(\theta-T_1)L}(\widetilde{U},\widetilde{V})
(h_{ \varepsilon,\delta})(\cdot,T_1) 
- e^{-(\theta-T_1)L}(\tilde{U}_1,\tilde{V}_1)]\Vert_{(L^2(\Omega))^2}
$$
$$
\le \kappa_6(\theta-T_1)^{-m}
\Vert(\widetilde{U},\widetilde{V})(h_{ \varepsilon,\delta})(\cdot,T_1)
- (\tilde{U}_1, \tilde{V}_1)\Vert_{(L^2(\Omega))^2}
\le \kappa_6(\theta-T_1)^{-m}( \varepsilon + \kappa_5\delta ).
$$
Extending $h_{ \varepsilon,\delta}(\cdot,t) = 0$ for $t > T_1$, we have
$$
e^{-(\theta-T_1)L}(\widetilde{U},\widetilde{V})(h_{ \varepsilon,\delta}) 
(\cdot,T_1) = (\widetilde{U},\widetilde{V})(h_{ \varepsilon,\delta})
(\cdot,\theta), 
$$
and so
\begin{equation}
\Vert L^m[(\tilde{U},\tilde{V})(h_{ \varepsilon,\delta})(\cdot,\theta)
- e^{-(\theta-T_1)L}(\tilde{U}_1,\tilde{V}_1)]\Vert_{(L^2(\Omega))^2}
\le \kappa_6(\theta-T_1)^{-m}( \varepsilon + \kappa_5\delta).
                                                    \label{(4.14)}
\end{equation}
Moreover as $(\tilde{U}_1, \tilde{V}_1) \in D(L^{m+1})$, we have
\begin{eqnarray*}
&&\Vert L^m[ e^{-(\theta-T_1)L}(\tilde{U}_1,\tilde{V}_1) 
- (\tilde{U}_1, \tilde{V}_1)]\Vert_{(L^2(\Omega))^2}
\le \kappa_6\Vert (e^{-(\theta-T_1)L} - I)L^m(\tilde{U}_1,
\tilde{V}_1)\Vert_{(L^2(\Omega))^2}\\
&\leq &\kappa_6\left\Vert \int^{\theta-T_1}_0 \frac{d}{d\eta}
(e^{-\eta L})L^m(\tilde{U}_1,\tilde{V}_1) d\eta 
\right\Vert_{(L^2(\Omega))^2}
\leq \kappa_6\left\Vert \int^{\theta-T_1}_0 e^{-\eta L}
L^{m+1}(\tilde{U}_1,\tilde{V}_1) d\eta \right\Vert
_{(L^2(\Omega))^2}\\
&\le &\kappa_6(\theta-T_1)\Vert L^{m+1}(\tilde{U}_1,\tilde{V}_1)
\Vert_{(L^2(\Omega))^2}.
\end{eqnarray*}
In terms of  (\ref{(4.14)}), we obtain
$$
\Vert L^m[(\tilde{U},\tilde{V})(h_{ \varepsilon,\delta})(\cdot,\theta)
- (\tilde{U}_1,\tilde{V}_1)] \Vert_{(L^2(\Omega))^2}
\le \kappa_6(\theta-T_1)^{-m}( \varepsilon + \kappa_5\delta )
$$
\begin{equation}
+ \kappa_6(\theta-T_1)\Vert L^{m+1}(\tilde{U}_1,\tilde{V}_1)\Vert
_{(L^2(\Omega))^2}.                              \label{(4.15)} 
\end{equation}
For any $ \varepsilon_1 > 0$ and $(\tilde{U}_1, \tilde{V}_1) \in 
D(L^{m+1})$, we choose $T_1 \in (0, \theta)$
such that
$$
\kappa_6(\theta-T_1)\Vert L^{m+1}(\tilde{U}_1,\tilde{V}_1)
\Vert_{(L^2(\Omega))^2} < \frac{ \varepsilon_1}{3}.
$$  
Then, with this $T_1$, we choose $ \varepsilon > 0$ such 
that 
$$
\kappa_6(\theta - T_1)^{-m} \varepsilon < \frac{ \varepsilon_1}{3}.
$$ 
Finally with this $h_{ \varepsilon}$, we choose $\delta > 0$ sufficiently small such that
$$
\kappa_6\delta\leq  \frac{ \varepsilon_1}{3}.
$$
Therefore  (\ref{(4.15)}) yields
\begin{equation}
\Vert L^m[(\tilde{U},\tilde{V})(h_{ \varepsilon,\delta})(\cdot,\theta)
- (\tilde{U}_1,\tilde{V}_1)]\Vert_{L^2(\Omega)^2}
 \leq  \varepsilon_1.                       \label{(4.16)}
\end{equation}
In terms of Lemma 4.3 and  (\ref{(4.1)}), by choosing $ \varepsilon > 0$ 
sufficiently small for $\inf_{x\in \Omega\setminus\omega_1}
\vert \tilde{U}_1(x)\vert$ and
$\inf_{x\in \Omega\setminus\omega_1}\vert \tilde{V}_1(x)\vert$,
the proof of Lemma 4.2 is completed.
\\
{\bf Second Step}
We will complete the proof of Theorem 4.1.
Let $h\in L^2(\omega_T)$ be chosen in Lemma 4.2.
We set
$$ 
u = U - \widetilde{U}, \quad v = V - \widetilde{V}.
$$
Then $(u,v)$ satisfies
$$
\partial_t u = \Delta u + au + bv + A\cdot\nabla u +
B\cdot \nabla v + f\tilde{V}, 
$$
$$
\partial_t v = \Delta v + cu + dv + C\cdot\nabla u +
D\cdot \nabla v + g\tilde{U} \quad \mbox{in $\Omega_T$},
$$
\begin{equation}
u=v=0 \qquad \mbox{on $\Sigma_T$}, \label{(4.17)}
\end{equation}
where 
$$
f = b - \widetilde{b}, \quad 
g = c - \widetilde{c}.
$$
We consider the time derivative of system  (\ref{(4.17)}).  Setting
$y=\partial_t u$ and $z=\partial_t v$, we obtain
$$
  \partial_t y= \Delta y+a(x)y+b(x)z+A\cdot\nabla y+B\cdot\nabla z
+ f\partial_t\widetilde{V} \quad \mbox{in $\Omega_T$},
$$
$$
\partial_t z= \Delta z+c(x)y+d(x)z+C\cdot\nabla y+D\cdot\nabla z
+ g\partial_t\widetilde{U}  \quad \mbox{in $\Omega_T$}, 
$$
\begin{equation}
y=z=0 \qquad \mbox{on $\Sigma_T$}.                     \label{(4.18)}
\end{equation}
Applying the Carleman estimate Theorem 2.2 to system  (\ref{(4.17)}) and
using $f=0$ in $\omega_1$, we have
$$
\int_{\Omega_T} (s\rho)^{-1}e^{-2s\eta_{\omega}}
(\vert \partial_ty\vert^2 + \vert \partial_tz\vert^2
+ \vert \Delta y\vert^2 + \vert \Delta z\vert^2
$$
$$
+ (s\rho)^2\vert\nabla y\vert^2 + (s\rho)^2\vert\nabla z\vert^2
+ (s\rho)^4\vert y\vert^2 + (s\rho)^4\vert z\vert^2)dxdt
$$
\begin{equation}
\le \kappa_7(s)\Vert y\Vert^2_{W^{2,1}_2(\omega_T)}
+ \kappa\int \hspace{-6.5pt} \int_{\Omega_T} e^{-2s \eta_{\omega}} 
(|f\partial_t\widetilde{V}|^2+|g\partial_t\widetilde{U}|^2) dxdt.
                                        \label{(4.19)}
\end{equation}
Furthermore, for large $s>0$, we can prove that
\begin{equation}
\int_{\Omega_T} \vert f(x)\vert^2 e^{-2s\eta_{\omega}(x,t)}dxdt
\le o(1)\int_{\Omega} \vert f(x)\vert^2 e^{-2s\eta_{\omega}
(x,\theta)} dx \quad \mbox{as $s \to \infty$}.     \label{(4.20)}
\end{equation}
In fact, we can prove similarly to \cite{IY:98}.  Recall that 
$T = 2\theta$.  Setting $\ell(t) = t(T-t)$, by (\ref{(2.2)}) we have 
$\frac{\partial(-\eta_{\omega})}{\partial t}(x,\theta) = 0$, $x \in \Omega$ and
$$
\frac{\partial^2(-\eta_{\omega})}{\partial t^2}(x,t)
= -\alpha_{\omega}(x)\frac{2\ell'(t)^2 - \ell(t)\ell''(t)}{\ell(t)^3}
$$
and
$$
\frac{\partial^3(-\eta_{\omega})}{\partial t^3}(x,t)
= -\alpha_{\omega}(x)\frac{6\ell'(t)(\ell(t)\ell''(t) - \ell'(t)^2)}
{\ell(t)^4}, \quad (x,t) \in \Omega_T.
$$
Therefore
$$
\frac{\partial^2(-\eta_{\omega})}{\partial t^2}(x,t)
\le -\frac{\kappa_8}{\ell(t)^3}, \quad (x,t) \in \Omega_T
$$
with a positive constant $\kappa_8$ and
$$
\frac{\partial^3(-\eta_{\omega})}{\partial t^3}(x,t) \ge 0, \quad 
0 \le t \le \theta, \thinspace x \in \Omega,
$$
$$
\frac{\partial^3(-\eta_{\omega})}{\partial t^3}(x,t) \le 0, \quad 
\theta \le t \le T, \thinspace x \in \Omega.
$$
Consequently by the mean value theorem, we can take $t_1$ such that
$t_1$ is between $t$ and $\theta$ and
$$
-\eta_{\omega}(x,t) = -\eta_{\omega}(x,\theta) 
+ \frac{1}{2}\frac{\partial^2(-\eta_{\omega})}{\partial t^2}(x,t)(t-\theta)^2
+ \frac{1}{6}\frac{\partial^3(-\eta_{\omega})}{\partial t^3}(x,t_1)(t-\theta)^3
$$
$$
\le -\eta_{\omega}(x,\theta) - \frac{\kappa_8}{2t^3(T-t)^3}(t-\theta)^2,
\quad (x,t) \in \Omega_T.
$$
Hence, noting that $\kappa_8 > 0$ and $-\frac{1}{t(T-t)} 
\le -\frac{4}{T^2}$, we obtain
$$
\int^T_0 e^{-2s\eta_{\omega}(x,t)} dt 
\le e^{-2s\eta_{\omega}(x,\theta)}
\int^T_0 \exp\left( -\frac{s\kappa_8}{t^3(T-t)^3}(t-\theta)^2\right) 
dt
$$
$$
\le e^{-2s\eta_{\omega}(x,\theta)}
\int^T_0 \exp\left( -\frac{s\kappa_9}{T^2}(t-\theta)^2\right) dt.
$$
The Lebesgue theorem yields .

We have
$$
\int_{\Omega} (\vert y(x,\theta)\vert^2 + \vert z(x,\theta)\vert^2)
e^{-2s\eta_{\omega}(x,\theta)} dx
= \int_{\Omega} \frac{\partial}{\partial t}
\int^{\theta}_0 (\vert y(x,t)\vert^2 + \vert z(x,t)\vert^2)
e^{-2s\eta_{\omega}(x,t)} dtdx
$$
$$
= \int_{\Omega} \int^{\theta}_0 
\{ 2s\alpha_{\omega}(x)\rho(t)^2(T-2t)
(\vert y(x,t)\vert^2 + \vert z(x,t)\vert^2)
+ 2(y\partial_ty + z\partial_tz)\} e^{-2s\eta_{\omega}(x,t)} dtdx
$$
$$
\le \kappa_{10}\int_{\Omega_T} 
\{ (s\rho)^3(\vert y(x,t)\vert^2 + \vert z(x,t)\vert^2)
+ (s\rho)^{-1}(\vert \partial_ty(x,t)\vert^2 + \vert \partial_tz(x,t)\vert^2)
\} e^{-2s\eta_{\omega}(x,t)} dtdx.
$$
At the last inequality, we used 
$$
\vert y\partial_ty\vert 
= \vert (s\rho)^{-\frac{1}{2}}\partial_ty(s\rho)^{\frac{1}{2}}y\vert
$$
$$
\le \frac{1}{2}(s\rho)^{-1}\vert \partial_ty\vert^2
+ \frac{1}{2}(s\rho)\vert y\vert^2
\le \frac{1}{2}(s\rho)^{-1}\vert \partial_ty\vert^2
+ \kappa'_{10}(s\rho)^3\vert y\vert^2.
$$
Hence, by  (\ref{(4.19)}) and  (\ref{(4.20)}), noting that $f=g=0$ in $\omega_1$, we have
\begin{equation}
\int_{\Omega} (\vert y(x,\theta)\vert^2 + \vert z(x,\theta)\vert^2)
e^{-2s\eta_{\omega}(x,\theta)} dx
\le \kappa_7(s)\Vert y\Vert^2_{W^{2,1}_2(\omega_T)}
+ o(1)\int_{\Omega\setminus\omega_1} 
(\vert f(x)\vert^2 + \vert g(x)\vert^2)
e^{-2s\eta_{\omega}(x,\theta)} dx                     \label{(4.21)}
\end{equation}
for all large $s > 0$.

On the other hand, since $u(\cdot,\theta) = v(\cdot,\theta) = 0$, we
have $y(x,\theta) = f(x)\tilde{V}(x,\theta)$ and
$z(x,\theta) = g(x)\tilde{U}(x,\theta)$ for $x \in \Omega$.
Therefore, by  (\ref{(4.6)}) and  (\ref{(4.21)}) we obtain
$$
\kappa_{11}\int_{\Omega} (\vert f(x)\vert^2 + \vert g(x)\vert^2)
e^{-2s\eta_{\omega}(x,\theta)} dx
\le \kappa_7(s)\Vert y\Vert^2_{W^{2,1}_2(\omega_T)}
+ o(1)\int_{\Omega\setminus\omega_1} 
(\vert f(x)\vert^2 + \vert g(x)\vert^2)
e^{-2s\eta_{\omega}(x,\theta)} dx
$$
as $s\longrightarrow \infty$.
Taking $s > 0$ large and fixing, we absorb the second term on the 
right hand side into the left hand side and the proof of
Theorem 4.1 is completed.
\section{Some  generalization and comments}
\label{sec: 5}
\subsection{Identification of all the coefficients}
%

\noindent
Indeed we can determine  all the coefficients of (\ref{syst-UV}). 
For it, we need repeats of measurements by choosing suitable interior 
controls.  We choose $m \in {\bf N}$ such that
$$
m > \frac{n}{4} + \frac{1}{2}.
$$
We recall that $(\tilde{U},\tilde{V}) = (\tilde{U}(h)(\cdot,\cdot),
\tilde{V}(h)(\cdot,\cdot))$ satisfies  (\ref{(4.3)}) and that 
$(U,V) = (U(h)(\cdot,\cdot), V(h)(\cdot,\cdot))$ satisfies  (\ref{(4.3)}) 
where $\tilde{a}, \tilde{b}, \tilde{c}, \tilde{d}, \tilde{A},
\tilde{B}, \tilde{C}, \tilde{D}$ are replaced by $a,b,c,d,A,B,C,D$
respectively.  Then, with $m > \frac{n}{4}+\frac{1}{2}$, 
under assumption  (\ref{(4.2)}) we 
can prove (\ref{(4.16)}).  Moreover, noting that $H^{2m}(\Omega)
\subset C^1(\overline{\Omega})$ by $m > \frac{n}{4}+\frac{1}{2}$,
we can see that for any $\varepsilon > 0$ and 
$(\tilde{U}_1, \tilde{V}_1) \in (H^{2m+2}_0(\Omega))^2$, there exists
$h \in L^2(\omega_T)$ such that
\begin{equation}
\Vert (\tilde{U}(h),\tilde{V}(h))(\cdot,\theta) 
- (\tilde{U}_1, \tilde{V}_1)\Vert_{(C^1(\overline{\Omega}))^2}
\le \varepsilon.                                          \label{(5.1)}
\end{equation}
Therefore
\\
{\bf Theorem 5.1}
{\it 
Let $\omega$, $a,\widetilde{a},b,\widetilde{b},c,\widetilde{c},d,
\widetilde{d},A,\widetilde{A},B,\widetilde{B},C,\widetilde{C},D,
\widetilde{D}$ satisfy Assumption 1.1 and 
$\widetilde{a},\widetilde{b},\widetilde{c},
\widetilde{d},\widetilde{A},\widetilde{B},\widetilde{C},
\widetilde{D} \in W^{2m,\infty}(\Omega)$.  
Let $\omega_1$ be a neighbourhood of $\partial \Omega$ such that
$\overline{\omega} \subset \omega_1$ and let the coefficients 
$(a,b,c,d,A,B,C,D)$ and $(\widetilde{a},\widetilde{b},\widetilde{c},
\widetilde{d},\widetilde{A},\widetilde{B},\widetilde{C},\widetilde{D})$ 
coincide in $\omega_1$. 
Then there exist $h_1, h_2, ..., h_{2n+2} \in L^2(\omega_T)$ 
such that
$$
\mbox{det}\thinspace
\left(\begin{array}{cccccccc}
\widetilde{U}(h_1) & \widetilde{V}(h_1)    & 0 & 0
&   \nabla \widetilde{U}(h_1)  & \nabla  \widetilde{V}(h_1) &0& 0 \\
0 & 0& \widetilde{U}(h_1) & \widetilde{V}(h_1) &0 & 0  
& \nabla \widetilde{U}(h_1)& \nabla  \widetilde{V}(h_1) \\
\widetilde{U}(h_2) & \widetilde{V}(h_2) &0  & 0
&   \nabla \widetilde{U}(h_2) & \nabla \widetilde{V}(h_2)  &0& 0 \\
0   & 0 &\widetilde{U}(h_2)  & \widetilde{V}(h_2)  & 0  &0
& \nabla \widetilde{U}(h_2) & \nabla \widetilde{V}(h_2) \\
\vdots & \vdots & \vdots & \vdots & \vdots &\vdots & \vdots & \vdots\\
\widetilde{U}(h_{2n+2}) & \widetilde{V}(h_{2n+2})    & 0 & 0
&   \nabla \widetilde{U}(h_{2n+2})  
& \nabla  \widetilde{V}(h_{2n+2}) &0& 0 \\
0 & 0& \widetilde{U}(h_{2n+2}) & \widetilde{V}(h_{2n+2}) &0 & 0  
& \nabla \widetilde{U}(h_{2n+2})& \nabla  \widetilde{V}(h_{2n+2})\\
\end{array} \right)
$$
\begin{equation}
\neq 0  \qquad x \in \overline{\Omega\setminus\omega_1}, \quad
t = \theta
                                           \label{(5.2)}
\end{equation}
and we choose a constant $\kappa > 0$ depending on
$M,m,\gamma, s,  \Omega,\omega,T$  and $h_1, ..., h_{2n+2}$ such that
\begin{eqnarray*}
&\Vert a - \widetilde{a}\Vert_{L^2(\Omega)}
+ \Vert b - \widetilde{b}\Vert_{L^2(\Omega)}
+ \Vert c - \widetilde{c}\Vert_{L^2(\Omega)}
+ \Vert d - \widetilde{d}\Vert_{L^2(\Omega)}\\
+ & \Vert A - \widetilde{A}\Vert_{(L^2(\Omega))^n}
+ \Vert B - \widetilde{B}\Vert_{(L^2(\Omega))^n}
+ \Vert C - \widetilde{C}\Vert_{(L^2(\Omega))^n}
+ \Vert D - \widetilde{D}\Vert_{(L^2(\Omega))^n}\\
\le &\kappa \sum_{j=1}^{2n+2} 
(\Vert \partial_t(U(h_j) - \tilde{U}(h_j))\Vert_{W_2^{2,1}(\omega_T)}
+ \Vert U(h_j) - \tilde{U}(h_j)\Vert_{W_2^{2,1}(\omega_T)})
\end{eqnarray*}
for all $(a,b,c,d,A,B,C,D)$ satisfying Assumption 1.1.}
\\  
{\bf Example for Theorem 5.1}:\\
Let $n=1$ and let $p_1, p_2, q_1, q_2, q_3$ be constants such that
$p_1q_2 - p_2q_1 \ne 0$ and $p_3(x_1)$, $q_4(x_1)$ satisfy
$(\partial_1p_3)(x_1) \ne 0$ and $\partial_1q_4(x_1) \ne
0$ for $x_1 \in \overline{\Omega\setminus \omega_1}$, and let $q_3$ be 
an arbitrarily smooth function.
Then for $x = x_1 \in \overline{\Omega\setminus\omega_1}$, we can 
verify that 
$$
\mbox{det}\;
\left(\begin{array}{cccccccc}
p_1 & q_1   & 0 & 0   &   \partial_1p_1 & \partial_1q_1   & 0 & 0 \\
0   & 0 & p_1 & q_1 &   0     &0     & \partial_1p_1&  \partial_1q_1\\
p_2 & q_2 &0 & 0   &   \partial_1p_2 & \partial_1q_2 &0& 0 \\
0   & 0 & p_2 & q_2 &   0  &0        & \partial_1p_2 & \partial_1q_2\\
\vdots & \vdots & \vdots & \vdots & \vdots &\vdots & \vdots & \vdots\\
p_4 & q_4 & 0 & 0   &   \partial_1p_4 & \partial_1q_4 & 0 & 0 \\
0   & 0 &p_4 &  q_4 &   0 &0         & \partial_1p_4& \partial_1q_4\\
\end{array}\right)(x)
= \vert \partial_1p_3(x_1)\vert^2\vert \partial_1q_4(x_1)\vert^2(p_1q_2-p_2q_1)^2
\neq 0.
$$
Therefore in (\ref{(5.1)}), we can choose $(\tilde{U}_1, \tilde{V}_1)
= (p_j, q_j)$, $1 \le j \le 4$ to construct $h_1, h_2, h_3, h_4$
satisfying (\ref{(5.2)}).

\subsection{Carleman estimate for a $3\times3$ reaction-diffusion system 
with one observation}
%

We consider now  a $3\times3$ reaction-diffusion system  
$$
\partial_t u(x,t)= \Delta u(x,t)+a_{11}(x,t)u(x,t)+a_{12}(x,t)v(x,t)
+ a_{13}(x,t)w(x,t) +f(x,t) \quad \mbox{in $\Omega_T$},
$$
$$
\partial_t v(x,t)= \Delta v(x,t)+a_{21}(x,t)u(x,t)+a_{22}(x,t)v(x,t)
+a_{23}(x,t)w(x,t)+g(x,t)  \quad \mbox{in $\Omega_T$}, 
$$
$$
\partial_t w(x,t)= \Delta w(x,t)+a_{31}(x,t)u(x,t)+a_{32}(x,t)v(x,t)
+ a_{33}(x,t)w(x,t) + h(x,t)\quad  \mbox{in $\Omega_T$}, 
$$
\begin{equation}
u=v=w=0 \qquad \mbox{on $\Sigma_T$}.             \label{(5.3)}
\end{equation}
We will assume\\
{\bf Assumption 5.2}\\
{\it 
(a) $(a_{ij})_{i,j=1,3} \in W^{2,\infty}(\Omega_T)$,
$\Vert a_{ij}\Vert_{W^{2,\infty}(\Omega_T)} \le M$.
\\
(b) $\omega $ of class $C^{2}$, $\partial \omega \cap \partial
\Omega =\gamma$ and $\left| \gamma \right| $ $\neq 0$.
\\
(c) $\left|(\nabla a_{12}-\frac{a_{12}}{a_{13}}\nabla a_{13})\cdot \nu
\right |\neq 0 \thinspace  \mbox{ on } \gamma\times (0,T)$.
\\
(d) $a_{12},a_{13}\in W^{3,\infty}(\omega_T)$,
$\Vert a_{12}\Vert_{W^{3,\infty}(\omega_T)}, 
\Vert a_{13}\Vert_{W^{3,\infty}(\omega_T)} \le M$.
\\
(e) $a_{13} \neq 0 \mbox { on } {\overline{\Omega_T}}$.
}

We show a Carleman estimate with extra data of one component.
\\  
{\bf Theorem 5.3}
{\it Under Assumption 5.2, there exist $\alpha_{\omega} \in 
C^2(\overline{\Omega})$ with $\alpha_{\omega} > 0$ on $\overline
{\Omega}$ and a constant $s_0>0$ which depends
on $T, M, \Omega, \omega, \tau$ and the $L^{\infty}(\Omega)$-norms
of $a_{ij}$, $1\le i,j \le 3$ such that we can choose positive 
constants $\kappa_1(s)$ and $\kappa$ satisfying:
$$
\int_{\Omega_T} (s\rho)^{-1}e^{-2s\eta_{\omega}}
(\vert \partial_t u\vert^2 + \vert \partial_t v\vert^2
+ \vert \partial_t w\vert^2 + \vert \Delta u\vert^2
+ \vert \Delta v\vert^2 + \vert \Delta w\vert^2
$$
$$
+ (s\rho)^2\vert \nabla u\vert^2 + (s\rho)^2\vert \nabla v\vert^2
+ (s\rho)^2\vert \nabla w\vert^2
+ (s\rho)^4 u^2 + (s\rho)^4 v^2 + (s\rho)^4 w^2) dxdt
$$
$$
\le \kappa_1(s)(\Vert u\Vert^2_{W^{4,2}_2(\omega_T)}
+ \Vert f\Vert^2_{W^{2,1}_2(\omega_T)}
+ \Vert g\Vert^2_{L^2(\omega_T)}
+ \Vert h\Vert^2_{L^2(\omega_T)})
$$
$$
+ \kappa\int_{\Omega_T} (\vert f^2\vert + \vert g\vert^2
+ \vert h\vert^2) e^{-2s\eta_{\omega}} dxdt
$$
for all $s \ge s_0$ and $(u,v,w)$ satisfying (\ref{(5.3)}).
Here we set 
$$
\eta_{\omega}(x,t) = \frac{\alpha_{\omega}(x)}{t(T-t)}.
$$
}
\\
{\bf Proof}
Setting $z=a_{12}v+a_{13}w$, we rewrite (\ref{(5.3)}) as
$$
\partial_t u= \Delta u+a_{11}u+z+f \quad \mbox{in $\Omega_T$},
$$
$$
\partial_t z= \Delta z+A\cdot \nabla z+az+eu+B\cdot \nabla v +bv+G \quad
\mbox{in $\Omega_T$}, 
$$
$$
\partial_t v= \Delta v+a_{21}u+ dv+cz+g \quad \mbox{in $\Omega_T$}, 
$$
\begin{equation}
u=v=z=0 \qquad \mbox{on $\Sigma_T$},            \label{(5.4)}
\end{equation}
where 
$$
A=-2 \frac{\nabla a_{13}}{a_{13}},\,
B=-2\nabla a_{12} +2\frac{a_{12}}{a_{13}}\nabla a_{13},
$$
$$ 
a = 2\frac{\vert \nabla a_{13}\vert^2}{a_{13}^2}
+ a_{33} + \frac{a_{12}a_{23}+\partial_ta_{13} - \Delta a_{13}}{a_{13}},
$$
$$
b = a_{12}a_{22} + a_{13}a_{32} + \partial_ta_{12} - \Delta a_{12}
+ 2\nabla a_{13}\cdot \nabla\left(\frac{a_{12}}{a_{13}}\right)
$$
$$
- \frac{a_{12}}{a_{13}}(a_{12}a_{23} + a_{13}a_{33} + \partial_ta_{13}
- \Delta a_{13}),
$$
$$
c = \frac{a_{23}}{a_{13}}, \quad 
d=a_{22} - \frac{a_{12}a_{23}}{a_{13}}, \quad
e=a_{21}a_{12}+a_{31}a_{13}
$$
and
$$
G=a_{12}g+a_{13}h.
$$ 

By \cite{FI}, \cite{I} and the proof of Theorem 2.2, we see that 
there exist a subdomain $\omega' \subset \omega$ and $\beta_{\omega'}
\in C^2(\overline{\Omega})$ with $\beta_{\omega'} > 0$ on $\overline
{\Omega}$ such that
$$
\int_{\Omega_T} (s\rho)^{-1}e^{-2s\tilde{\eta}_{\omega'}}
(\vert \partial_t u\vert^2 + \vert \partial_t z\vert^2
+ \vert \Delta u\vert^2
+ \vert \Delta z\vert^2 
$$
$$
+ (s\rho)^2\vert \nabla u\vert^2 + (s\rho)^2\vert \nabla z\vert^2
+ (s\rho)^4 u^2 + (s\rho)^4 z^2) dxdt
$$
$$
\le \kappa_2\int_{\Omega_T} (z^2 + \vert eu\vert^2) 
e^{-2s\tilde{\eta}_{\omega'}}dxdt
+ \kappa_2\int_{\Omega_T} (f^2 + G^2)e^{-2s\tilde{\eta}_{\omega'}}dxdt
$$
\begin{equation}
+ \kappa_2\int_{\Omega_T} \vert B\cdot \nabla b + bv\vert^2
e^{-2s\tilde{\eta}_{\omega'}}dxdt
+ \kappa_2\int_{\omega'_T} (s\rho)^{-1}(u^2 + z^2)
e^{-2s\tilde{\eta}_{\omega'}}dxdt                    \label{(5.5)}
\end{equation}
and
$$
\int_{\Omega_T} (s\rho)^{-1}e^{-2s\tilde{\eta}_{\omega'}}
(\vert \partial_t z\vert^2 + \vert \partial_t v\vert^2
+ \vert \Delta z\vert^2 + \vert \Delta v\vert^2 
$$
$$
+ (s\rho)^2\vert \nabla z\vert^2 + (s\rho)^2\vert \nabla v \vert^2
+ (s\rho)^4 z^2 + (s\rho)^4 v^2) dxdt
$$
\begin{equation}
\le \kappa_3(s)(\Vert z\Vert^2_{W^{2,1}_2(\omega'_T)}
+ \Vert eu + G\Vert^2_{L^2(\omega'_T)})
+ \kappa_2\int_{\Omega_T} (\vert eu+G\vert^2 + \vert a_{21}u+g\vert^2)
e^{-2s\tilde{\eta}_{\omega'}}dxdt              \label{(5.6)}
\end{equation}
for all $s \ge s_0$, where we set $\tilde{\eta}_{\omega'}(x,t) = 
\frac{\beta_{\omega'}(x)}{t(T-t)}$.
Here (\ref{(5.5)}) is obtained by applying the Carleman estimate in \cite{FI}
or \cite{I} to the first and the second equations in (\ref{(5.4)}), while
(\ref{(5.6)}) is seen by applying Theorem 2.2 to the second and the third
equations in (\ref{(5.4)}) and noting Assumption 5.2 (c).  We further notice
that the weight function $\tilde{\eta}_{\omega'}$ can be taken the 
same, which can be seen from the proof of Theorem 2.2.  By (\ref{(5.5)}) and 
(\ref{(5.6)}), in terms of Assumption 5.2 (a), (d) and (e), we have
$$
\int_{\Omega_T} (s\rho)^{-1}e^{-2s\tilde{\eta}_{\omega'}}
(\vert \partial_t u\vert^2 + \vert \partial_t z\vert^2 
+ \vert \partial_t v\vert^2 
+ \vert \Delta u\vert^2 + \vert \Delta z\vert^2 
+ \vert \Delta v\vert^2
$$
$$
+ (s\rho)^2\vert \nabla u\vert^2 + (s\rho)^2\vert \nabla z\vert^2
+ (s\rho)^2\vert \nabla v\vert^2 
+ (s\rho)^4 u^2 + (s\rho)^4 z^2 + (s\rho)^4v^2) dxdt
$$
$$
\le \kappa_2\int_{\Omega_T} (z^2 + u^2 + v^2 + \vert \nabla v\vert^2)
e^{-2s\tilde{\eta}_{\omega'}}dxdt
+ \kappa_2\int_{\Omega_T} (f^2 + g^2 + h^2)
e^{-2s\tilde{\eta}_{\omega'}}dxdt
$$
\begin{equation}+ \kappa_3(s)(\Vert z\Vert^2_{W^{2,1}_2(\omega'_T)}
+ \Vert u\Vert^2_{L^2(\omega'_T)} + \Vert g\Vert^2
_{L^2(\omega'_T)} + \Vert h\Vert^2_{L^2(\omega'_T)})    \label{(5.7)}
\end{equation}
for all large $s > 0$.  
We can absorb the first terms on the right hand side into the left
hand side by choosing $s > 0$ large, and we use
$z = \partial_t u - \Delta u - a_{11}u - f$ by the first equation in (\ref{(5.4)}),
so that the proof of Theorem 5.3 is completed.

The approximate controllability is a direct consequence of Theorem 5.3.
That is, we consider 
$$
\partial_t u= \Delta u+a_{11}(x) u+ a_{21}(x)v+a_{31}(x) w 
+\chi_{\omega}f \quad \mbox{in $\Omega_T$},
$$
$$
\partial_t v= \Delta v+a_{12}(x)u+a_{22}(x) v+a_{32}(x) w \quad
\mbox{in $\Omega_T$},
$$
$$
\partial_t w= \Delta w+a_{13}(x)u+a_{23}(x)v+a_{33}(x)w \quad
\mbox{in $\Omega_T$}, 
$$
$$
u=v=w=0 \qquad \mbox{on $\Sigma_T$},
$$
\begin{equation}
u(\cdot,0)=u_0,\quad v(\cdot,0)=v_0,\quad w(\cdot,0)=w_0 \quad
\mbox{in $\Omega$}.                   \label{(5.8)}
\end{equation}
Here we assume that all the coefficients are independent of 
$t$.

Then
\\
{\bf Theorem 5.4}
{\it
Under Assumption 5.2, for all $\varepsilon >0,\,T>0$, 
$(u_0,v_0,w_0)\in (L^2(\Omega))^3$ and
$(u_1,v_1,w_1)\in (L^2(\Omega))^3$, there exists $f\in L^2(\omega_T)$ 
such that the corresponding solution of (\ref{(5.8)}) satisfies 
$$\Vert(u,v,w)(\cdot,T)-(u_1,v_1,w_1)\Vert_{(L^2(\Omega))^3 }
\leq \varepsilon.
$$
}

Similarly to section 4, we can apply the Carleman estimate of Theorem 5.3
for determining the nine coefficients $a_{ij}$, $1\le i,j \le 3$ by
suitably repeated observations of only one component $u$ and we 
will here omit further details.

{\bf Acknowledgements.}  Most part of this paper was written 
during the stays of the fourth named author in 2007 at Universit\'{e}
de Provence and he thanks the invitation.\\
The fourth author was partially supported by l'Agence Nationale de la
Recherche under grant ANR JC07\_183284. 


\end{document}